\documentclass[10pt,thmsa]{article}
\usepackage{amsfonts, amsmath, latexsym, amssymb, amsthm, amscd}
\usepackage{graphicx}

\newtheorem{theorem}{Theorem}
\newtheorem{corollary}{Corollary}
\newtheorem{proposition}{Proposition}
\newtheorem{lemma}{Lemma}
\newtheorem{definition}{Definition}
\newtheorem{remark}{Remark}
\newcommand{\p}{\Bbb{P}}

\newcommand{\ind}{\mbox{\rm 1\hspace{-0.04in}I}}

\newcommand{\R}{\mbox{\rm I\hspace{-0.02in}R}}

\newcommand{\eqdef}{\stackrel{\mbox{\tiny$($def$)$}}{=}}

\def\QED{\hfill\vrule height 1.5ex width 1.4ex depth -.1ex \vskip20pt}

\newcommand{\ud}{\mathrm{d}}
\begin{document}
\hspace*{-0.5in}{\footnotesize This version March 6, 2008.}
\vspace*{0.9in}
\begin{center}
{\Large On the Lamperti stable  processes}\\
\vspace*{0.4in} {\large M.E. Caballero\footnote{$^{, 3}$ Instituto
de Matem\'aticas, Universidad Nacional Autonoma de M\'exico,
M\'exico D.F  C.P. 04510. $^1$E-mail: marie@matem.unam.mx
$^{3}$E-mail: garmendia@matem.unam.mx}, J.C. Pardo\footnote{
Department of Mathematical Science, University of Bath. {\sc Bath}
BA2 7AY. {\sc United Kingdom}, E-mail:jcpm20@bath.ac.uk} and J.L.
P\'erez$^3$\vspace*{0.2in}}\\

\end{center}
\vspace{0.2in}
\begin{abstract}
\par We consider a new family  of  $\R^d$-valued L\'{e}vy processes  that we call Lamperti stable. One of the advantages of this class is that the law of many related functionals can be computed explicitly. Here, in the one dimensional case we provide an explicit form for the characteristic exponent and other several useful properties of the class (for some particular cases see \cite{cc}, \cite{ckp}, \cite{kp} and \cite{pp}). This family of processes shares many properties with the tempered stable  and the layered stable  processes, defined in Rosi\'nski \cite{ro} and Houdr\'e and Kawai \cite{hok} respectively, for instance their short and long time behavior. We also find a series representation which is used for sample path simulation, illustrated in the case $d=1$. Finally we provide many examples, some of which appear in recent literature.

\noindent {\sc Key words and phrases:} Lamperti stable distributions and processes, stable processes, L\'evy processes.\\

\noindent MSC 2000 subject classification: 60E07, 60G51, 60G52.

\end{abstract}
\section{Introduction.}
In recent years the interest in having more
accurate models in various domains of applied probability  has lead
to an increasing attention paid to some special classes of L\'{e}vy
processes related to the stable law, for example: the tempered
stable  and the layered stable processes introduced by Rosi\'{n}ski
in \cite{ro} and Houdr\'{e} and Kawai in \cite{hok}, respectively.
Both families of processes have nice structural and analytical properties, such as
combining in short time the behavior of stable processes and in
long time the behavior of a Brownian motion. They also have a
series representation which may be used for sample paths simulation.

 Lamperti \cite{la} and more
recently, Caballero and Chaumont \cite{cc} studied four families of
L\'evy processes  which are related to the stable subordinator and
some conditioned stable processes via the Lamperti representation of
positive self-similar Markov processes. Those studies are the
starting point of our work. Recall that positive self-similar Markov
processes,  $(X,\p_x)$, $x>0$, are strong Markov process with
c\`adl\`ag paths, which fulfill a scaling property, i.e.~there
exists a constant $\alpha > 0$ such that for any $b>0$:
\[
\mbox{\it The law of $\;(bX_{b^{-\alpha}t},\,t\ge0)$ under $\p_x$ is
$\p_{bx}$.}
\]
We shall refer to these processes as pssMp. According to Lamperti
\cite{la}, any pssMp up to its first hitting time of 0 may be
expressed as the exponential of a L\'evy process, time changed by
the inverse of its exponential functional. Reciprocally, any L\'evy process
$\xi$  can be expressed as the logarithm of a time changed pssMp $X$.
In this paper we refer to this Lamperti transformation as $LT_1$ and the details can be seen in \cite{la}.

One of the examples treated by Lamperti in \cite{la} is the case when
 $(X,\mathbf{P}_x)$ is a stable subordinator of index
$\alpha\in(0,1)$ starting from $x>0$. Lamperti in  \cite{la}, describes the
characteristics of the associated L\'evy process $\mathfrak{s}=(\mathfrak{s}_t,
t\ge 0)$ which is again a subordinator, with no drift and with
L\'evy measure given by
\[
\eta(\ud x)=\frac{c_+e^x}{(e^x-1)^{\alpha+1}}\ud x,  \qquad x > 0.
\]
The three cases of pssMp studied in \cite{cc} are related to some conditioned stable processes. The
first one is the stable L\'evy processes killed when it first exits
from the positive half-line, here denoted by $(X^*,\mathbf{P}_x)$.
The second class corresponds to  that of stable processes
conditioned to stay positive (see for instance \cite{ch, do}),
denoted by $(X^\uparrow,\mathbf{P}_x)$. Finally, the third class of
pssMp   is that of stable processes conditioned to hit 0
continuously, denoted by $(X^\downarrow,\mathbf{P}_x)$. The
corresponding L\'evy processes under the $LT_1$ transformation are
denoted by $\xi^*,\xi^\uparrow$ and $\xi^\downarrow$, respectively.
These three classes of L\'evy processes have no gaussian component and their L\'evy measure are of the type
\[
\pi(\ud x)=\left(\frac{c_+e^{bx}}{(e^x-1)^{\alpha+1}}\ind_{\{x>0\}}+
\frac{c_-e^{ bx}}{(1-e^x)^{\alpha+1}}\ind_{\{x<0\}}\right)\, \ud x,
\]
where $c_+, c_-$ are the constants of the L\'evy measure of the original stable process and $b$ is a positive parameter. We recall that for $\xi^*$ the constant $b$ equals 1 and moreover it has finite lifetime if $c_->0$ and,  when $c_-=0$, it has infinite lifetime and drifts to $-\infty$. For $\xi^{\uparrow}$ the constant $b$ is equal to $\alpha\rho+1$, where $\rho=\mathbf{P}_0(X_1<0)$. It has infinite lifetime and  drifts to $\infty$.  Finally, for the processes $\xi^{\uparrow}$ the constant $b$ is $\alpha\rho$. It has infinite lifetime and  drifts to $-\infty$. We remark that such processes have linear coefficients in their respective characteristic exponent that we denote by  $a^*, a^{\uparrow}$ and $a^{\downarrow}$. Such constants are computed
explicitly in \cite{cc} in terms of $\alpha$, $\rho$, $c_-$ and
$c_+$. Actually it was proved recently in \cite{ckp}, that the
process $\xi^\downarrow$ corresponds to $\xi^\uparrow$ conditioned
to drift to $-\infty$ (or equivalently $\xi^\uparrow$ is
$\xi^\downarrow$ conditioned to drift
to $+\infty$).

Finally, motivated by self-similar continuous state branching
processes with immigration, Patie has recently studied in \cite{pp}
the family of L\'evy process with no positive jumps with L\'evy
measure
\[
\eta^*(\ud x)=\frac{c_- e^{(\alpha+\vartheta)
x}}{(1-e^x)^{\alpha+1}}\ud x, \qquad x < 0,
\]
where $\vartheta>-\alpha$.

These L\'evy processes  have the advantage that the law of many
functionals  can be computed explicitly, for example: the first exit
time from a finite interval or semi-finite interval, overshoots
distributions and exponential functionals  (see for instance
Caballero and Chaumont \cite{cc}, Chaumont et al. \cite{ckp},
Kyprianou and Pardo \cite{kp} and Patie \cite{pp}). We also
emphasize that in some cases the Wiener-Hopf factors can be computed
and scale functions in the spectrally one sided case can be
obtained. Since many tractable mathematical expressions can be
computed, this class seems to be an useful tool for applications and
rich enough to be of particular interest.

In this work, we investigate a generalization  of the L\'evy
processes mentioned above and we will refer to them as Lamperti
stable processes.
We also study these processes in higher
dimensions. We will see that
this class has nice structural and analytical properties close to those for tempered stable and layered stable processes.

In section 2 we begin by studying the Lamperti stable  distributions, which are multivariate infinitely divisible
distributions with no Gaussian component and whose L\'{e}vy measure is characterized by a
triplet $(\alpha,f,\sigma)$, more precisely an
index $\alpha\in (0,2)$, a function $f$, and a finite measure $\sigma$, both defined on the unit sphere in $\mathbb{R}^{d}$.
In particular, the radial component of
any of these L\'evy measures is asymptotically equivalent to that of an
stable distribution, with index $\alpha$, near zero and has exponential decay at
infinity. These distributions have a
density with respect to the Lebesgue measure, have finite moments of all orders and  exponential finite
moments of some order. In the one dimensional case, the density is
$C^{\infty}$. In some particular cases, we also prove that these
distributions are  self-decomposable.

In section 3, we formally introduce the Lamperti stable processes
and study their properties with emphasis in the one dimensional
case, where we obtain an explicit closed form for the characteristic exponent.
Motivated by the works of Rosi\'nski \cite{ro} and Houdr\'e
and Kawai \cite{hok}, we  prove in section 4, 5 and 6 that Lamperti
stable processes in a short time  look like a stable process while
in a large time scale  they look like a Brownian motion,  that they
are  absolute continuous  with respect to  its short time limiting
stable process and they admit a series representation that allows
simulations of their paths, respectively.

In the last section, we
study some related processes:  the Ornstein-Uhlenbeck processes
whose limiting distribution is a Lamperti stable law and the L\'evy
processes with no positive jumps whose descending ladder height
process is a Lamperti stable subordinator. Finally we illustrate  with several examples the presence of
Lamperti stable distributions in recent literature.

\section{Lamperti stable distributions.}
In this section, we  define  Lamperti stable distributions  on $\R^d$ and establish some of their basic properties. According to Theorem
14.3 in Sato \cite{sa}, the L\'evy measure $\Pi$ of a stable
distribution with index $\alpha$ on $\R^d$ in polar coordinates is
of the form
\[
\Pi(dr,d\xi)=r^{-(\alpha+1)}\ud r\sigma(\ud \xi)
\]
where $\alpha\in(0,2)$ and $\sigma$ is a finite measure on $S^{d-1}$, the unit sphere on $\R^d$.  The measure $\sigma$ is uniquely determined by $\Pi$. Conversely, for any non-zero finite measure $\sigma$ on $S^{d-1}$ and for any $\alpha\in (0,2)$ we can   define an stable distribution with L\'evy measure defined as above.

Motivated by the form of the L\'evy measure of the processes
mentioned in the introduction  and the previous discussion, we define a
new family of  infinitely divisible distributions that we call
Lamperti stable.
\begin{definition}
Let $\mu$ be an infinitely divisible probability measure on $\R^{d}$
without Gaussian component. Then, $\mu$ is called Lamperti stable
if its L\'{e}vy measure on $\R^d_0:=\R^{d}\setminus\{0\}$ is given
by
\begin{equation}\label{e1}
\nu_{\sigma}^{\alpha,f}(B)=\int_{S^{d-1}}\sigma(\ud\xi)\int_{0}^{\infty}\ind_{B}(r\xi)e^{rf(\xi)
}(e^{r}-1)^{-(\alpha+1)}dr,\quad\quad B\in \mathcal{B}(\R_{0}^{d}),
\end{equation}
where $\alpha\in(0,2)$, $\sigma$ is non-zero finite measure on $S^{d-1}$, and
$f:S^{d-1}\to\mathbb{R}$ is a measurable function such that $\gamma:=\sup_{\xi\in S^{d-1}}f(\xi)<\alpha+1$.
\end{definition}
Note that  $\nu_{\sigma}^{\alpha,f}$ is indeed a L\'{e}vy
measure on $\R_{0}^{d}$.  To see this, we need to verify that
\begin{equation}\label{levmea}
\int_{\scriptsize
\R_{0}^{d}}(1\wedge\|x\|^{2})\nu_{\sigma}^{\alpha,f}(\ud
x)=\int_{S^{d-1}}\sigma(\ud\xi) \int_{0}^{\infty}(1\wedge
r^{2})e^{rf(\xi)}(e^{r}-1)^{-(\alpha+1)}\ud r<\infty.\notag
\end{equation}
On the one hand, since
$e^{rf(\xi)}(e^{r}-1)^{-(\alpha+1)}\sim r^{-(\alpha+1)}$ as $ r\to 0$, \footnote{We say that $f\sim g$ as $x\rightarrow x_{0}$ if for $x_{0}\in \mathbb{R}^{d}$ $\lim_{x \to x_{0}}f(x)/g(x)=1$.} we have  that
\[
\int_{S^{d-1}}\sigma(\ud\xi) \int_{0}^{1}
r^{2}e^{rf(\xi)}(e^{r}-1)^{-(\alpha+1)}\ud r<\infty.
\]
One the other hand, from elementary calculations we deduce
\[
\int_{S^{d-1}}\sigma(\ud\xi) \int_{1}^{\infty}
e^{rf(\xi)}(e^{r}-1)^{-(\alpha+1)}\ud r\leq \frac{\sigma(S^{d-1})}{(1-e^{-1})^{\alpha+1}}\int_1^\infty e^{-r(\alpha+1-\gamma)}\ud r,
\]
where $\gamma:=\sup_{\xi\in S^{d-1}}f(\xi)$. Since $\gamma<\alpha+1$, the above integral is finite and therefore $\nu_{\sigma}^{\alpha,f}$ is a L\'evy measure.

In the one dimensional case $f$ takes only two possible values, since $S^0=\{-1,1\}$. In the sequel, we denote these two values by $f(1):=\beta$ and $f(-1):=\rho$ . From the definition of
the measure $\sigma$, we have $\sigma(\{1\})=c_+$ and
$\sigma(\{-1\})=c_-$. Therefore each distribution associated  to the
L\'evy processes mentioned in the introduction belongs to the class
of Lamperti stable distribution. Following the notation of the
introduction, for:
\begin{itemize}
\item the subordinator $\mathfrak{s}$,  $\beta=1$ and $c_-=0$,
\item  the process $\xi^*$, $\beta=1$ and $\rho=\alpha$,
\item  the process $\xi^{\uparrow}$, $\beta=\alpha\rho+1$ and $\rho=\alpha(1-\rho)$,
\item  the process $\xi^{\downarrow}$, $\beta=\alpha\rho$ and $\rho=\alpha(1-\rho)+1$,
\item  the class of L\'evy processes with no positive jumps considered by Patie \cite{pp}, $\rho=1-\vartheta$ and $c_+=0$.
\end{itemize}
Note that Lamperti stable distributions satisfy the
divergence condition, i.e.
\[
\int_{0}^{\infty}e^{rf(\xi)}(e^{r}-1)^{-(\alpha+1)}dr=\infty \qquad
\textrm{for any}\qquad \xi\in S^{d-1}.
\]
Thus from  Theorem 27.10 in \cite{sa}, we deduce that they are
absolutely continuous with respect to the Lebesgue measure. Also note
that Lamperti stable distributions do not belong in general to the
class of
 tempered $\alpha$-stable distributions (see \cite{ro}), nor to the class of layered stable distributions (see \cite{hok}
 for their definition). This follows respectively from the facts that the function
\begin{equation}
q(r,\xi)=\frac{e^{rf(\xi)}}{(e^{r}-1)^{\alpha+1}}r^{1+\alpha}\notag
\end{equation}
is not completely monotone for any fixed $\xi\in S^{d-1}$, and
that for each $\xi\in S^{d-1}$
\begin{equation}
\frac{q(r,\xi)}{r^{\alpha+1}}\sim
e^{-(\alpha+1-f(\xi))r}\qquad\hbox{as }\quad r\to\infty.\notag
\end{equation}
The above remarks motivate  us  to see Lamperti
stable distributions as a  particular class  and study some of their properties.
\begin{proposition}
Let $\mu$ be a Lamperti stable distribution with L\'{e}vy measure
$\nu_{\sigma}^{\alpha,f}$ given by (\ref{e1}), then
\begin{equation}
\int_{\scriptsize\R^{d}}\|x\|^{p}\mu(dx)<\infty\quad \textrm{ for
all }\quad p>0.\notag
\end{equation}
Let $\gamma:=\sup_{\xi\in S^{d-1}} f(\xi)$. If $\zeta
<\alpha+1-\gamma$, then \[
\int_{\scriptsize\R^{d}}e^{\zeta\|x\|}\mu(dx)<\infty,
\]
In particular, for $\kappa<\alpha+1$ and if $f\equiv \kappa$, we
have
\[
\int_{\scriptsize\R^{d}}e^{\zeta\|x\|}\mu(dx)<\infty\qquad \textrm{
if and only if }\qquad \zeta<\alpha+1-\kappa.
\]
\end{proposition}
\noindent {\it Proof:} According to Theorem 25.3 in Sato \cite{sa},
we need to verify  that the restriction of $\nu^{\alpha,f}_{\sigma}$
to the set $\{x\in\R_{0}^{d}:\|x\|>1\}$ has the corresponding moment
properties. To this end, we consider
\begin{align*}
\int_{\{\|x\|>1\}}\|x\|^{p}\nu^{\alpha,f}_{\sigma}(dx)&=\int_{S^{d-1}}\sigma(d\xi)\int_{1}^{\infty}r^{p}e^{f(\xi)
r}(e^{r}-1)^{-(\alpha+1)}dr\notag\\
&\leq
\sigma(S^{d-1})(1-e^{-1})^{-(\alpha+1)}\int_{1}^{\infty}r^{p}e^{r(\gamma-(\alpha+1))}dr,\notag
\end{align*}
which is finite since  $\gamma<\alpha+1$.\\
Now, we turn our attention to the exponential moments. Consider
\begin{align}\label{a:8}
\int_{\{\|x\|>1\}}e^{\zeta\|x\|}\nu^{\alpha,f}_{\sigma}(dx)&=\int_{S^{d-1}}\sigma(d\xi)\int_{1}^{\infty}e^{r\zeta}e^{f(\xi)
r}(e^{r}-1)^{-(\alpha+1)}dr\notag\\
&\leq
\sigma(S^{d-1})(1-e^{-1})^{-(\alpha+1)}\int_{1}^{\infty}e^{r(\zeta+\gamma-(\alpha+1))}dr,\notag
\end{align}
which is also finite for $\zeta<\alpha+1-\gamma$.  Next, we suppose that
$f\equiv\kappa$.  The former arguments imply that for
$\zeta<\alpha+1-\kappa$, the Lamperti stable distribution $\mu$ has
a finite exponential moment of order $\zeta$. In a similar way, it
is clear that
\begin{equation}
\int_{\{\|x\|>1\}}e^{\zeta\|x\|}\nu^{\alpha,\kappa}_{\sigma}(dx)\geq\sigma(S^{d-1})\int_{1}^{\infty}e^{r(\zeta+\kappa-(\alpha+1))}dr.\notag
\end{equation}
This implies that
$\displaystyle\int_{\{\|x\|>1\}}e^{\zeta\|x\|}\nu^{\alpha,\kappa}_{\sigma}(dx)$
is finite if and only if $\zeta<\alpha+1-\kappa$.\QED

Our next result shows that Lamperti stable distributions are selfdecomposable in some cases, and that they belong to
the Jurek class, i.e. the class of infinitely divisible distributions for which the L\'{e}vy measure $\nu$ takes the following form
\begin{equation}
\nu(B)=\int_{S^{d-1}}\sigma(\ud\xi)\int_{0}^{\infty}\ind_{B}(r\xi)l(\xi,r)dr,\quad\quad \text{for}\quad B\in \mathcal{B}(\R_{0}^{d}),\notag
\end{equation}
where $l(\xi,r)$ is measurable $\xi\in S^{d-1}$, and decreasing in $r\in(0,\infty)$.
\begin{proposition} Let $\mu$ be an Lamperti-stable distribution
on $\R^{d}$ with L\'{e}vy measure $\nu^{\alpha,f}_{\sigma}$ given by
(\ref{e1}), then $\mu$ belong to the Jurek class. Moreover, $\mu$ is
selfdecomposable if and only if $f(\xi)\leq \alpha+1/2$, for all
$\xi\in S^{d-1}$ and $\alpha\in (0,2)$.
\end{proposition}
\noindent {\it Proof:} We first prove that  Lamperti stable
distributions belong to the Jurek class. To do so, it is
enough to show that the function
\[
\ell(\xi,r)=\frac{e^{f(\xi) r}}{(e^{r}-1)^{\alpha+1}},
\]
is measurable in $\xi\in S^{d-1}$ and decreasing in $r>0$. The
measurability of $\ell$ in $\xi\in S^{d-1}$ is clear from its
definition. To prove that $\ell$ is decreasing in $r>0$, we fix $\xi\in S^{d-1}$ and consider the derivative of $\ell_1(\cdot)=\ell(\xi,\cdot)$,  i.e.
\[
\ell_1'(r)=\frac{e^{f(\xi)
r}}{(e^{r}-1)^{\alpha+2}}\Big(e^r(f(\xi)-\alpha-1)-f(\xi)\Big).
\]
Hence $\ell_1'(r)<0$ for $r>0$, since $f(\xi)\leq \alpha+1$. This implies that $\ell_1(r)$ is decreasing for $r>0$ and prove the statement.

Now, we prove the second part of the Proposition. From
Theorem 15.10 in Sato \cite{sa} we  need to verify that  the function
\begin{equation}
k(\xi,r)=\frac{e^{f(\xi) r}}{(e^{r}-1)^{\alpha+1}}r,\notag
\end{equation}
is measurable in $\xi\in S^{d-1}$ and decreasing for $r>0$, if and only if
 $f(\xi)\leq\alpha+1/2$ for all
$\xi\in S^{d-1}$ and $\alpha\in(0,2)$. The measurability of $k$ in $\xi\in S^{d-1}$ is clear from its
definition. Fix $\xi\in S^{d-1}$, and assume that
$f(\xi)\leq \alpha+1/2$. Consider the derivative of $h(r)=k(\xi,r)$,
\begin{equation}
h'(r)=\frac{e^{f(\xi)
r}}{(e^{r}-1)^{\alpha+2}}\Big(e^{r}\big(1+f(\xi)r-(\alpha+1)r\big)-f(\xi)
r-1\Big)\notag.
\end{equation}
Therefore, the function $h$ is  decreasing for $r>0$ if
\begin{equation}
g(r)=e^{r}(1+(f(\xi)-(\alpha+1))r)-f(\xi) r-1,\notag
\end{equation}
is negative for $r>0$. Again, straightforward calculations give us
that
\begin{align}
g'(r)&=e^{r}\Big(f(\xi)-\alpha+f(\xi)r-(\alpha+1)r\Big)-f(\xi)\quad\textrm{and}\notag\\
g''(r)&=e^{r}\Big(2f(\xi)-2\alpha-1+f(\xi)r-(\alpha+1)r\Big).\notag
\end{align}
It is not difficult to see that for $f(\xi)\leq \alpha+1/2$, we have
that $g''(r)<0$ for $r>0$, which implies that $g'(r)$ is decreasing
for $r>0$. Now, since  $g'(0)=-\alpha<0$, we deduce   that the
function $g'(r)<0$ for $r>0$, which implies that $g$ is decreasing
for $r>0$. Next, since $g(0)=0$, we have that the function $g$ is
negative for
$r>0$, which proves that $h$ is decreasing for $r>0$ assuming that $f(\xi)\leq \alpha+1/2$, for all $\xi\in S^{d-1}$.\\
Finally, it just remains to prove that if $f(\xi)>\alpha+1/2$, for
some $\xi\in S^{d-1}$ and $\alpha\in (0,2)$, there exists a
$r_{0}>0$ such that $g(r_{0})>0$. This follows  from similar
arguments used  above (take for instance
$f(\xi)\equiv\alpha+b+1/2$, for $b>0$, and $\alpha$ close to 0).\QED We finish this
section with some properties of Lamperti stable distributions
defined on $\R$. The first of which says in particular that the
density of any Lamperti stable distribution belongs to $C^{\infty}$.
\begin{proposition}
Let $\mu$ be a Lamperti stable distribution on $\R$, then $\mu$ has
a $C^{\infty}$ density and all the derivatives of the density tend
to $0$ as $|x|$ tends to $\infty$.
\end{proposition}
\noindent{\it Proof:} Recall that the function $f$
takes two values,  $\beta=f(1)$ and $\rho=f(-1)$
as usual.  According to Proposition
28.3 in Sato \cite{sa} and the form of the L\'evy measure of $\mu$,
it is enough to prove that
\begin{equation}\label{cond0}
g(r)=\int_{0}^r x^2 \frac{e^{\beta x}}{(e^x-1)^{\alpha+1}}\ud x,
\quad\textrm{verifies that}\quad
\liminf_{r\to 0}\frac{g(r)}{r^{2-a}}>0,
\end{equation}
for some $a\in(0,2)$. But this is immediate because for $r$
sufficiently small, we have
\[
\int_{0}^r x^2 \frac{e^{\beta x}}{(e^x-1)^{\alpha+1}}\ud x\geq K
\int_0^r\frac{x^2}{x^{\alpha+1}}\ud x=K r^{2-\alpha},
\]
where $K>0$. In particular when  $a=\alpha$, the condition in
(\ref{cond0}) is satisfied and the statement follows.\QED Before we
state the last result of this section, we  recall the definition of
a particular class of distributions  which is important in risk
theory (see for instance \cite{ekm} and \cite{kkm}).
\begin{definition}[Class $\mathcal{L}^{(\delta)}$] Take a parameter $\delta\geq0$. We
shall say that a distribution function $G$ on $[0,\infty)$ with tail
$\overline{G}:=1-G$ belongs to class $\mathcal{L}^{(\delta)}$ if
$\overline{G}(x)>0$ for each $x\geq0$ and
\begin{equation}
\lim_{u\to\infty}\frac{\overline{G}(u-x)}{\overline{G}(u)}=e^{\delta
x}\quad\hbox{for each $x\in\mathbb{R}$}.\notag
\end{equation}
The tail of any (L\'evy or other) measure, finite and non-zero on
$(x_0,\infty)$ for some $x_0>0$, can be renormalised to be the tail
of a distribution function and by extension, then is said to be in
$\mathcal{L}^{(\delta)}$, if the associated distribution function is in $\mathcal{L}^{(\delta)}$.
\end{definition}
Now, we will prove that the tail of the L\'evy measure of  any
Lamperti stable distributions defined in $\R$ belongs to the class
$\mathcal{L}^{(\alpha+1-\beta)}$, where $\beta=f(1)$ as usual.
\begin{proposition}
Let $\mu$ be a Lamperti-stable distribution on $\R$, then the tail
of its L\'evy measure belongs to the class
$\mathcal{L}^{(\alpha+1-\beta)}$. In particular when $\mu$ is
defined on $\R_+$,  we have that $\mu$ belongs to the class
$\mathcal{L}^{(\alpha+1-\beta)}$.
\end{proposition}
\noindent{\it Proof:} First, we define
\begin{equation}
\nu(u)=\frac{1}{K}\int_{1}^{u}\frac{e^{\beta
r}}{(e^{r}-1)^{\alpha+1}}\ud r, \qquad u\geq 1,\notag
\end{equation}
where $K=\displaystyle\int_{1}^{\infty}\frac{e^{\beta
r}}{(e^{r}-1)^{\alpha+1}}\ud r$. Note that $\nu$ corresponds to the
distribution function associated to the tail of the L\'{e}vy measure
of a
Lamperti stable distribution.\\
From elementary calculations, we get
\begin{align}
\frac{\overline{\nu}(u-x)}{\overline{\nu}(u)}&=\int_{u-x}^{\infty}\frac{e^{\beta
r}}{(e^r-1)^{\alpha+1}}\ud r\left(\int_{u}^{\infty}\frac{e^{\beta
r}}{(e^{r}-1)^{\alpha+1}}\ud r\right)^{-1}\notag\\
&\leq
(\alpha+1-\beta)e^{u(\alpha+1-\beta)}\int_{u-x}^{\infty}\frac{e^{\beta
r}}{(e^{r}-1)^{\alpha+1}}\ud r\notag\\
&=e^{(\alpha+1-\beta)x}\Big(1-e^{-(u-x)}\Big)^{-\alpha-1}.\notag
\end{align}
Similarly
\begin{align}
\frac{\overline{\nu}(u-x)}{\overline{\nu}(u)} &\geq
(1-e^{-u})^{\alpha+1}(\alpha+1-\beta)e^{u(\alpha+1-\beta)}\int_{u-x}^{\infty}e^{(\beta-(\alpha+1))
r}\ud r\notag\\
&=e^{(\alpha+1-\beta)x}(1-e^{-u})^{\alpha+1}.\notag
\end{align}
Therefore taking $u$ large enough, we deduce that
$\nu\in\mathcal{L}^{(\alpha+1-\beta)}$. The  case when $\mu$ is defined in $\R_+$ follows from Proposition 3.4 in
Kyprianou et al. \cite{kkm}. \QED
\section{Lamperti stable L\'{e}vy processes.}
Here, we introduce the class of L\'{e}vy processes which is
associated to Lamperti stable distributions. We also discuss,
specially in the one dimensional case, a number of coarse and fine
properties of their paths which are of particular interest for
applications.
\begin{definition}
A L\'{e}vy process without gaussian component, and linear term $\theta$, is called Lamperti
stable with characteristics $(\alpha,f,\sigma,\theta)$ if its L\'{e}vy
measure is given by (\ref{e1}).
\end{definition}
In  the sequel, we denote the Lamperti stable L\'{e}vy process with
characteristics $(\alpha,f,\sigma,\theta)$  by $X^{L}=(X^{L}_{t},t\geq0)$.
Its characteristic exponent is defined by $E[\exp(i\langle y,X_{t}^{L}\rangle)]=\exp(-t\Psi(y))$ for $t\geq0$, $y\in \R^d$ where
\begin{equation}\label{lk}
\Psi(y)=i\langle
y,\theta\rangle+\int_{\scriptsize\R^{d}_{0}}\Big(1-e^{i\langle
y,x\rangle}+i\langle
y,x\rangle\ind_{\{\|x\|<1\}}\Big)\nu^{\alpha,f}_{\sigma}(\ud
x),
\end{equation}
the measure $\nu^{\alpha,f}_{\sigma}$ has the form given in
(\ref{e1}) and $\theta\in \R^{d}$.

The first  property in study is the $p$-th  variation  of Lamperti
stable processes. In particular, we prove that their $p$-th
variation is similar to that of stable processes.
\begin{proposition}
 Let $X^{L}$ be a Lamperti stable process with characteristics $(\alpha,f,\sigma,\theta)$.
\begin{itemize}
\item[i)] If $\alpha\in(1,2)$, the process $X^{L}$ is a.s. of finite $p$-th  variation in every finite interval if and only if $p\in(\alpha,2)$.
\item[ii)] The process $X^{L}$ is a.s. of finite  variation in every finite interval if
and only if $\alpha\in(0,1)$.
\end{itemize}
\end{proposition}
\noindent{\it Proof:}  $(i)$ From Theorem III in Bretagnolle
\cite{bre}, we have that for $p\in(1,2)$, the process $X^{L}$ is a.s.
of finite $p$-th  variation on every finite interval if and only if
\[
\int_{\{\|x\|\leq1\}}\|x\|^{p}\nu^{\alpha,f}_{\sigma}(\ud x)<\infty.
\]
Recall that $\gamma:=\sup_{\xi\in S^{d-1}} f(\xi)$. From the form of
the L\'evy measure $\nu^{\alpha,f}_{\sigma}$ and some elementary
calculations, we have
\begin{equation}\label{eq1}
\int_{\{\|x\|\leq1\}}\|x\|^{p}\nu^{\alpha,f}_{\sigma}(\ud x)\leq
\sigma(S^{d-1})e^{\gamma}\int_{0}^1\frac{r^p}{(e^r-1)^{\alpha+1}}\ud
r\leq \sigma(S^{d-1})e^{\gamma}\int_{0}^1r^{p-(\alpha+1)}\ud r.
\end{equation}
On the other hand, we have
\begin{equation}\label{eq2}
\begin{split}
\int_{\{\|x\|\leq1\}}&\|x\|^{p}\nu^{\alpha,f}_{\sigma}(\ud x)\geq \sigma\Big(\{\xi\in S^{d-1}:f(\xi)\geq 0\}\Big)\int_{0}^1\frac{r^p}{(e^r-1)^{\alpha+1}}\ud r\\
&+\int_{S^{d-1}}\ind_{\{f(\xi)< 0\}}e^{f(\xi)}\sigma(\ud x)\int_{0}^1\frac{r^p }{(e^r-1)^{\alpha+1}}\ud r\\
&\geq K\bigg(\sigma\Big(\{\xi\in S^{d-1}:f(\xi)\geq
0\}\Big)+\int_{S^{d-1}}\ind_{\{f(\xi)< 0\}}e^{f(\xi)}\sigma(\ud
x)\bigg)\int_0^1r^{p-(\alpha+1)}\ud r,
\end{split}
\end{equation}
for some $K>0$. Therefore
$X^{L}$ is of finite $p$-th variation on every finite interval if and only if $p>\alpha$.\\
The proof of part $(ii)$ is very similar. According to Theorem 3 of
Gikhman and Skorokhod \cite{gsk}, it  is enough to prove that
\[
\int_{\{\|x\|\leq1\}}\|x\|\nu^{\alpha,\gamma}_{\sigma}(dx)<\infty,
\]
if and only if $\alpha\in(0,1)$. But this follows from  (\ref{eq1})
and (\ref{eq2}) taking $p=1$, which concludes the proof.\QED Recall
that the characteristic exponent of a L\'evy process has a simpler
expression when its sample paths  have a.s. finite variation in
every finite interval. In this case, we have that the characteristic exponent in (\ref{lk}) takes
the form
\[
\Psi_{L}(y)=-i\langle \mathbf{d},y\rangle+\int_{\scriptsize
\R^d_0}\Big(1-e^{i\langle y,x\rangle}\Big)\nu_\sigma^{\alpha,f}(\ud x),
\]
where $\displaystyle\mathbf{d}=-\theta-\int_{\{\|x\|\leq1\}}x\nu_\sigma^{\alpha,f}(\ud x)$ and, in this case, $\mathbf{d}$ is known as the drift
coefficient.

In the rest of this section we work with real valued processes. We are now interested in two
important properties of L\'evy processes: creeping and the
regularity of $0$. Define for each $x\geq 0$, the first
passage time
\[
\tau^+_x=\inf\Big\{t>0: X^{L}_t>x\Big\},
\]
with the convention $\inf \emptyset=\infty$. We say that the
Lamperti stable process $X^{L}$ {\it creeps upwards} if for all
$x\geq 0$, $\mathbf{P}_0(X^{L}_{\tau^+_x}=x)>0$.  If $-X^{L}$ creeps
upwards, we  say that $X^{L}$ {\it creeps downwards}.
\begin{proposition}
Let $X^{L}$ be a Lamperti stable process   with
characteristics $(\alpha, f, \sigma, \theta)$.
\begin{itemize}
\item[i)]  If $\alpha\in(0,1)$ and  $\mathbf{d}>0$, the process $X^{L}$ creeps upwards.
\item[ii)] If $\alpha\in (1,2)$ the process $X^{L}$ does not creeps upwards.
\end{itemize}
\end{proposition}
\noindent{\it Proof:} $(i)$ The statement follows directly  from  part $(i)$ of Theorem 8 in \cite{kl}.\\
 $(ii)$ From Theorem 8 part $(iii)$ in \cite{kl}, all we have to
 prove is:
\begin{equation}\label{ii1}
\int_{0}^1\frac{x\nu^{\alpha,f}_{\sigma}\big([x,\infty)\big)}{H(x)}\ud
x=\infty, \quad\textrm{where} \quad
H(x)=\int_{-x}^0\int_{-1}^{y}\nu^{\alpha,f}_{\sigma}\big((-\infty,u]\big)\ud
u \ud y.
\end{equation}
To this end, we study the behavior of the integrand  in (\ref{ii1})
near $0$. We first note that
\begin{equation}\label{ii2}
u^{\alpha}\int_{-\infty}^{u}\frac{e^{-\rho
x}}{(e^{-x}-1)^{\alpha+1}}\ud r\sim 1\qquad\textrm{as }\quad
u\uparrow 0,
\end{equation}
which follows from  the fact that
$e^{-x\rho}(e^{-x}-1)^{-(\alpha+1)}\sim x^{-(\alpha+1)}$ when $x$ increases to $0$. Thus, it is not difficult
to deduce that
\[
x^{\alpha-2}H(x)\sim \frac{1}{(2-\alpha)(\alpha-1)\alpha}\qquad
\textrm{as }\quad x\downarrow 0.
\]
On the other hand, similar arguments as those used in (\ref{ii2})
give us
\begin{equation}\label{ii3}
x^{\alpha}\nu^{\alpha,f}_{\sigma}\big([x,\infty)\big)\sim\frac{1}{\alpha}\qquad
\textrm{as }\quad x\downarrow 0.
\end{equation}
Hence,
\[
\frac{x\nu^{\alpha,f}_{\sigma}\big([x,\infty)\big)}{H(x)}\sim
(2-\alpha)(\alpha-1)\frac{1}{x}\qquad \textrm{as }\quad x\downarrow
0,
\]
which implies (\ref{ii1}). The proof is now complete.\QED
\begin{proposition}
For  a Lamperti stable process $X^{L}$  with characteristics
$(\alpha,f,\sigma, \theta)$, the point $0$ is regular for $(0,\infty)$ if one of these two conditions hold:
\begin{itemize}
\item[i)] $\alpha\in [1,2)$.\\
\item[ii)] $\alpha\in(0,1)$, and $\mathbf{d}\geq0$.
\end{itemize}
\end{proposition}
\noindent\textit{Proof}:  Here, we apply Theorem 11 in \cite{kl}.
When $\alpha\in[1,2)$, the process $X^{L}$ has unbounded variation
(see Proposition 5 $(ii)$). Hence from part $(i)$ of Theorem 11 in
\cite{kl} we deduce that $0$ is regular for
$(0,\infty)$.\\
Now, we study the case where $\alpha\in(0,1)$.  In this case, the
process $X^{L}$ has bounded variation and according to part $(ii)$
and $(iii)$ of Theorem 11 in \cite{kl}, the point $0$ is regular for
$(0,\infty)$ if the drift coefficient $\mathbf{d}>0$ or if $\mathbf{d}=0$
and the following condition holds
\begin{equation}\label{cond1}
\int_{0}^{1}\frac{x\nu^{\alpha,f}_\sigma(\ud x)}{H(x)}=\infty,\quad
\textrm{where}\quad H(x)=\int_{0}^x
\nu^{\alpha,f}_\sigma(-\infty,-y)\ud y.
\end{equation}
Thus, we only need to verify  (\ref{cond1}). In order to do so, we
recall (\ref{ii3})
\[
y^{\alpha}\nu^{\alpha,f}_{\sigma}((-\infty,y])\sim\frac{1}{\alpha}\qquad\textrm{as}\quad
y\downarrow 0,
\]
which implies that
\[
x^{\alpha-1}H(x)\sim\frac{1}{\alpha(1-\alpha)}\qquad\textrm{as}\quad
x\downarrow 0.
\]
We observe then,  that
\[
\frac{x^{2}e^{\beta
x}(e^{x}-1)^{-(\alpha+1)}}{H(x)}\sim\alpha(1-\alpha),\qquad\textrm{as}\quad
x\downarrow 0,
\]
which implies (\ref{cond1}) and  the proof is now complete.\QED
Our
next result deals with the computation of  the characteristic  exponents of Lamperti stable processes. Denote by
 \[
 (z)_{\alpha}=\frac{\Gamma(z+\alpha)}{\Gamma(z)},\qquad \textrm{for} \quad z\in\mathbb{C},
 \]
 which is known as the Pochhammer symbol.
\begin{theorem}
Let $X^{L}$ be a Lamperti stable process with characteristics $(\alpha,f,\sigma,\theta)$.\\
\begin{itemize}
\item[i)]If $\alpha\in(0,1)\cup(1,2)$, the characteristic exponent of $X^{L}$ is given by
\begin{align}
\Psi_{L}(\lambda)&=i\lambda\tilde{\theta}-c_{+}\Gamma(-\alpha)\left((-i\lambda+1-\beta)_{\alpha}-(1-\beta)_{\alpha}\right)\notag\\
&-c_{-}\Gamma(-\alpha)\left((i\lambda+1-\rho)_{\alpha}-(1-\rho)_{\alpha}\right),\quad\quad\lambda\in\mathbb{R}.\notag
\end{align}
\item[ii)]If $\alpha=1$, the characteristic exponent of $X^{L}$ is given by
\begin{align}
\Psi_{L}(\lambda)&=i\lambda\tilde{\theta}-c_{+}\bigg((-i\lambda+1-\beta)\psi(-i\lambda+2-\beta)-(1-\beta)\psi(2-\beta)\bigg)\notag\\
&-c_{-}\bigg((i\lambda+1-\rho)\psi(i\lambda+2-\rho)-(1-\rho)\psi(2-\rho)\bigg),\quad\quad\lambda\in\mathbb{R}.\notag
\end{align}
\end{itemize}
Where $\tilde{\theta}$ is given by
\begin{equation}\label{exp:40}
\tilde{\theta}=\\
 \left\{
              \begin{array}{ll}
                -\mathbf{d}& \hbox{if $\alpha\in(0,1)$,}\\
                \displaystyle\theta-\left(c_{+}\tilde{a}_{\beta}-c_{-}\tilde{b}_{\rho}+(c_{+}-c_{-})(1-\mathcal{C})\right) & \hbox{if $\alpha=1$,}\\
               \displaystyle\theta-\left(c_{+}\tilde{a}_{\beta}-c_{-}\tilde{b}_{\rho}+\frac{c_{+}-c_{-}}{\alpha-1}\right) & \hbox{if $\alpha\in(1,2)$,}
              \end{array}
\right.
\end{equation}
$\tilde{a}_{\beta},\tilde{b}_{\rho}$ are given in (\ref{exp:30}), (\ref{exp:31}), respectively; $\mathcal{C}$ is the Euler constant, and $\psi$ is the Digamma function.
\end{theorem}
\noindent {\it Proof}: $i)$ First we will consider the case where $\alpha\in(0,1)$. Without loss of generality, we assume that $\mathbf{d}=0$. Since $\alpha\in(0,1)$, we know that the characteristic exponent of $X^{L}$ is given by
\begin{equation}\label{exp:1}
\Psi_{L}(\lambda)=-\left(c_{+}\int_{0}^{\infty}(e^{i\lambda x}-1)\frac{e^{\beta x}}{(e^{x}-1)^{\alpha+1}}dx+c_{-}
\int_{-\infty}^{0}(e^{i\lambda x}-1)\frac{e^{-\rho x}}{(e^{-x}-1)^{\alpha+1}}dx\right).
\end{equation}
We compute each of these integrals which we call $I_{1}$ and $I_{2}$ respectively. Since all the computations involved are valid for all $\lambda\in\mathbb{R}$, we center our attention in the variable $\beta$. In order to compute $I_{1}$ explicitly we will define in the set $U=\{z\in\mathbb{C}:\Re(z)<\alpha+1\}$, the following function $F:U\rightarrow\mathbb{C}$, given by
\begin{align}\label{exp:2}
F(z)&:=\int_{0}^{\infty}(e^{i\lambda x}-1)\frac{e^{z x}}{(e^{x}-1)^{\alpha+1}}dx=\int_{0}^{\infty}(e^{i\lambda x}-1)\frac{e^{-z_{1} x}}{(1-e^{-x})^{\alpha+1}}dx\notag\\
&=\int_{0}^{1}(u^{-i\lambda}-1)u^{z_{1}-1}(1-u)^{-(\alpha+1)}du,
\end{align}
where $z_{1}=\alpha+1-z$ and $\Re(z_{1})>0$. Then by making an integration by parts in the last integral of (\ref{exp:2}) we obtain for $\Re(z_{1})>1$
\begin{align}\label{exp:3}
\int_{0}^{1}(u^{-i\lambda}-1)u^{z_{1}-1}(1-u)^{-(\alpha+1)}du&=\frac{(-i\lambda-z_{1}+1)}{\alpha}\int_{0}^{1}u^{-i\lambda+z_{1}-2}(1-u)^{-\alpha}du\notag\\
&+\frac{z_{1}-1}{\alpha}\int_{0}^{1}u^{z_{1}-2}(1-u)^{-\alpha}du.
\end{align}
Now recalling the integral representation for the Beta function, (see \cite{Le}), we have for $\Re(a),\Re(b)>0$
\begin{equation}\label{exp:21}
B(a,b)=\int_{0}^{1}u^{a-1}(1-u)^{b-1}du=\frac{\Gamma(a)\Gamma(b)}{\Gamma(a+b)}
\end{equation}
we can express (\ref{exp:3}), in the following form:
\begin{align}
\int_{0}^{1}(u^{-i\lambda}-1)u^{z_{1}-1}(1-u)^{-(\alpha+1)}du&=\frac{-(i\lambda+z_{1}-1)}{\alpha}\frac{\Gamma(-i\lambda+z_{1}-1)\Gamma(1-\alpha)}{\Gamma(-i\lambda+z_{1}-\alpha)}\notag\\
&+\frac{(z_{1}-1)}{\alpha}\frac{\Gamma(z_{1}-1)\Gamma(1-\alpha)}{\Gamma(z_{1}-\alpha)}\notag,
\end{align}
finally by the recurrence relation for the Gamma function, $\Gamma(x+1)=x\Gamma(x)$, and the fact that $z_{1}=\alpha+1-z$, we obtain
\begin{equation}\label{exp:4}
F(z)=\Gamma(-\alpha)\left(\frac{\Gamma(-i\lambda+\alpha+1-z)}{\Gamma(-i\lambda+1-z)}-\frac{\Gamma(\alpha+1-z)}{\Gamma(1-z)}\right),
\end{equation}
for $\Re(z_{1})>1$, i.e. $\Re(z)<\alpha$. So we have the desired result for $\beta<\alpha$. In order to obtain it for $\beta\in[\alpha,\alpha+1)$ we do the following:
The equality (\ref{exp:4}) is valid in particular in $D_{\alpha}=\{z\in\mathbb{C}:\|z\|<\alpha\}$, in order to extend it to the case where $\|z\|<\alpha+1$, we will prove first that both sides of the equality in (\ref{exp:4}) are analytic functions in the disk $D_{\alpha+1}=\{z\in\mathbb{C}:\|z\|<\alpha+1\}$.\\
First we take $F$, and then using a series expansion we have
\begin{align}\label{exp:5}
F(z)=\int_{0}^{\infty}(e^{i\lambda x}-1)\frac{e^{z x}}{(e^{x}-1)^{\alpha+1}}dx=\int_{0}^{\infty}\sum_{n=0}^{\infty}(e^{i\lambda x}-1)\frac{(z x)^{n}}{n!}(e^{x}-1)^{-(\alpha+1)}dx.
\end{align}
Now consider the following
\begin{align}\label{exp:22}
\int_{0}^{\infty}&\sum_{n=0}^{\infty}\left\|(e^{i\lambda x}-1)\frac{(z x)^{n}}{n!}(e^{x}-1)^{-(\alpha+1)}\right\|dx\notag\\
&\leq\int_{0}^{\infty}\sum_{n=0}^{\infty}(|\lambda|x)\frac{(\|z\| x)^{n}}{n!}(e^{x}-1)^{-(\alpha+1)}dx\notag\\
&=\int_{0}^{\infty}(|\lambda|x)e^{\|z\|x}(e^{x}-1)^{-(\alpha+1)}dx,
\end{align}
which is finite when $\|z\|<\alpha+1$, therefore we can apply Fubini's Theorem in (\ref{exp:5}) and obtain
\begin{align}\label{exp:23}
F(z)=\sum_{n=0}^{\infty}\frac{z^{n}}{n!}\int_{0}^{\infty}(e^{i\lambda x}-1)\frac{x^{n}}{(e^{x}-1)^{\alpha+1}}dx=\sum_{n=0}^{\infty}a_{n}z^{n},
\end{align}
for $z\in W$, where
\begin{equation}
a_{n}=\frac{1}{n!}\int_{0}^{\infty}(e^{i\lambda x}-1)\frac{x^{n}}{(e^{x}-1)^{\alpha+1}}dx.\notag
\end{equation}
which implies that $F$ is analytic in $D_{\alpha+1}$.\\
Since for $\|z\|<\alpha+1$ we have that $\Re(-i\lambda+\alpha+1-z)>0$, and $\Re(\alpha+1-z)>0$, therefore the function
$G:U\rightarrow\mathbb{C}$, given by
\begin{equation}
G(z)=\Gamma(-\alpha)\left(\frac{\Gamma(-i\lambda+\alpha+1-z)}{\Gamma(-i\lambda+1-z)}-\frac{\Gamma(\alpha+1-z)}{\Gamma(1-z)}\right),\notag
\end{equation}
is analytic en $D_{\alpha+1}$. Since $F$ and $G$ are analytic in $D_{\alpha+1}$, and $F\equiv G$ in $D_{\alpha}$, we conclude that $F\equiv G$ in $D_{\alpha+1}$, which implies that
\begin{equation}
I_{1}=F(\beta)=\Gamma(-\alpha)\left(\frac{\Gamma(-i\lambda+\alpha+1-\beta)}{\Gamma(-i\lambda+1-\beta)}-\frac{\Gamma(\alpha+1-\beta)}{\Gamma(1-\beta)}\right),\notag
\end{equation}
for all $\beta<\alpha+1$.\\
Now we compute the second integral in the right-hand side of (\ref{exp:1})
\begin{equation}
I_{2}=\int_{-\infty}^{0}(e^{i\lambda x}-1)\frac{e^{-\rho x}}{(e^{-x}-1)^{\alpha+1}}dx=\int_{0}^{1}(u^{i\lambda}-1)u^{\rho_{1}-1}(1-u)^{-(\alpha+1)}du,\notag
\end{equation}
where $\rho_{1}=\alpha+1-\rho$, hence following the same arguments used in the computation of $I_{1}$, we get
\begin{equation}
I_{2}=\Gamma(-\alpha)\left(\frac{\Gamma(i\lambda+\alpha+1-\rho)}{\Gamma(i\lambda+1-\rho)}-\frac{\Gamma(\alpha+1-\rho)}{\Gamma(1-\rho)}\right),\notag
\end{equation}
for all $\rho<\alpha+1$. Therefore from the form of $I_{1}$ and $I_{2}$ we get
\begin{align}
\Psi_{L}(\lambda)&=-c_{+}\Gamma(-\alpha)\left((-i\lambda+1-\beta)_{\alpha}-(1-\beta)_{\alpha}\right)\notag\\
&-c_{-}\Gamma(-\alpha)\left((i\lambda+1-\rho)_{\alpha}-(1-\rho)_{\alpha}\right).\notag
\end{align}
for all $\beta, \rho<\alpha+1$.\\
Now we consider the case where $\alpha\in(1,2)$. As in the case where $\alpha\in(0,1)$, we assume that $\theta=0$. Since $\alpha\in(1,2)$, the characteristic exponent of $X^{L}$ is given by
\begin{align}\label{exp:6}
\Psi_{L}(\lambda)&=-\bigg(c_{+}\int_{0}^{\infty}\left(e^{i\lambda x}-1-i\lambda x\ind_{\{x<1\}}\right)\frac{e^{\beta x}}{(e^{x}-1)^{\alpha+1}}dx\notag\\
&+c_{-}\int_{-\infty}^{0}\left(e^{i\lambda x}-1-i\lambda x\ind_{\{x>-1\}}\right)\frac{e^{-\rho x}}{(e^{-x}-1)^{\alpha+1}}dx\bigg),
\end{align}
We call $I_{1}$ and $I_{2}$ respectively the integrals in (\ref{exp:6}). To study $I_{1}$ to do that we define the function $G:U\rightarrow\mathbb{C}$, given by
\begin{align}\label{exp:11}
G(z)&:=\int_{0}^{\infty}\left(e^{i\lambda x}-1-i\lambda x\ind_{\{x<1\}}\right)\frac{e^{z x}}{(e^{x}-1)^{\alpha+1}}dx\notag\\
&=\int_{0}^{\infty}\left(e^{i\lambda x}-1-i\lambda x\ind_{\{x<1\}}\right)\frac{e^{-z_{1} x}}{(1-e^{-x})^{\alpha+1}}dx\notag\\
&=\int_{0}^{\infty}\frac{(e^{i\lambda x}-1)e^{-z_{1}x}-i\lambda(1-e^{-x})e^{-x}}{(1-e^{-x})^{\alpha+1}}dx+i\lambda\int_{0}^{1}\frac{x(e^{-x}-e^{-z_{1}x})}{(1-e^{-x})^{\alpha+1}}dx\notag\\
&+i\lambda \int_{0}^{1}e^{-x}\frac{1-x-e^{-x}}{(1-e^{-x})^{\alpha+1}}dx+i\lambda\int_{1}^{\infty}\frac{e^{-x}}{(1-e^{-x})^{\alpha}}dx\notag\\
&=i\lambda\tilde{a}+i\lambda I(z)+\int_{0}^{1}\frac{(u^{-i\lambda}-1)u^{z_{1}-1}+i\lambda(u-1)}{(1-u)^{\alpha+1}}du,
\end{align}
where $z_{1}=\alpha+1-z$, $\Re(z_{1})>0$,
\begin{equation}
\tilde{a}=\int_{0}^{1}e^{-x}\frac{1-x-e^{-x}}{(1-e^{-x})^{\alpha+1}}dx+\int_{1}^{\infty}\frac{e^{-x}}{(1-e^{-x})^{\alpha}}dx,\notag
\end{equation}
and $I:U\rightarrow\mathbb{C}$ is defined by
\begin{equation}\label{exp:12}
I(z):=\int_{0}^{1}\frac{xe^{-x}(1-e^{-(\alpha-z)x})}{(1-e^{-x})^{\alpha+1}}dx.
\end{equation}
We consider the last integral in (\ref{exp:11}), an integration by parts gives us for $\Re(z_{1})>2$
\begin{align}\label{exp:7}
\int_{0}^{1}\frac{(u^{-i\lambda}-1)u^{z_{1}-1}+i\lambda(u-1)}{(1-u)^{\alpha+1}}du&=\frac{i\lambda}{\alpha}-\frac{(z_{1}-1)}{\alpha}\int_{0}^{1}(u^{-i\lambda}-1)u^{z_{1}-2}(1-u)^{-\alpha}du\notag\\
&+\frac{i\lambda}{\alpha}\int_{0}^{1}(u^{-i\lambda+z_{1}-2}-1)(1-u)^{-\alpha}du.
\end{align}
We now compute the integrals in the right-hand side of (\ref{exp:7}), by making an integration by parts, recalling the integral form of the Beta function (\ref{exp:21}), and using the recursion formula for the Gamma function we obtain for the first integral the following
\begin{align}\label{exp:8}
\int_{0}^{1}&(u^{-i\lambda}-1)u^{z_{1}-2}(1-u)^{-\alpha}du\notag\\
&=\frac{(z_{1}-2)}{1-\alpha}\int_{0}^{1}(u^{-i\lambda}-1)u^{z_{1}-3}(1-u)^{1-\alpha}du-\frac{i\lambda}{1-\alpha}\int_{0}^{1}u^{-i\lambda+z_{1}-2}(1-u)^{1-\alpha}du\notag\\
&=\frac{(-i\lambda+z_{1}-2)}{1-\alpha}\int_{0}^{1}u^{-i\lambda+z_{1}-3}(1-u)^{1-\alpha}du-\frac{(z_{1}-2)}{1-\alpha}\int_{0}^{1}u^{z_{1}-3}(1-u)^{1-\alpha}du\notag\\
&=\Gamma(1-\alpha)\left(\frac{\Gamma(-i\lambda+z_{1})}{(-i\lambda+z_{1}-1)\Gamma(-i\lambda+z_{1}-\alpha)}-\frac{\Gamma(z_{1})}{(z_{1}-1)\Gamma(z_{1}-\alpha)}\right),
\end{align}
and for the second integral in the right-hand side of (\ref{exp:7})
\begin{align}\label{exp:9}
\int_{0}^{1}(u^{-i\lambda+z_{1}-2}-1)(1-u)^{-\alpha}du&=\frac{1}{\alpha-1}+\frac{(-i\lambda+z_{1}-2)}{1-\alpha}\int_{0}^{1}u^{-i\lambda+z_{1}-3}(1-u)^{1-\alpha}du\notag\\
&=\frac{1}{\alpha-1}+\frac{(-i\lambda+z_{1}-2)}{1-\alpha}\frac{\Gamma(-i\lambda+z_{1}-2)\Gamma(2-\alpha)}{\Gamma(-i\lambda+z_{1}-\alpha)}\notag\\
&=\frac{1}{\alpha-1}+\frac{\Gamma(-i\lambda+z_{1})}{(-i\lambda+z_{1}-1)}\frac{\Gamma(1-\alpha)}{\Gamma(-i\lambda+z_{1}-\alpha)},\notag\\
\end{align}
so using (\ref{exp:8}) and (\ref{exp:9}), in (\ref{exp:7}) and recalling that $z_{1}=\alpha+1-z$, we get
\begin{align}\label{exp:10}
\int_{0}^{1}\frac{(u^{-i\lambda}-1)u^{z_{1}-1}+i\lambda(u-1)}{(1-u)^{\alpha+1}}du=\frac{i\lambda}{\alpha-1}+\Gamma(-\alpha)\left(\frac{\Gamma(-i\lambda+\alpha+1-z)}{\Gamma(-i\lambda+1-z)}-\frac{\Gamma(\alpha+1-z)}{\Gamma(1-z)}\right),\notag
\end{align}
if $\Re(z)<\alpha-1$.\\
So if we consider the function $P:U\rightarrow\mathbb{C}$ defined by
\begin{equation}\label{exp:14}
P(z)=\frac{i\lambda}{\alpha-1}+\Gamma(-\alpha)\left(\frac{\Gamma(-i\lambda+\alpha+1-z)}{\Gamma(-i\lambda+1-z)}-\frac{\Gamma(\alpha+1-z)}{\Gamma(1-z)}\right),
\end{equation}
we have that
\begin{equation}
G(z)=i\lambda\tilde{a}+i\lambda I(z)+P(z)\notag,
\end{equation}
for $\Re(z)<\alpha-1$, in particular the equality holds in the set $D_{\alpha-1}=\{z\in\mathbb{C}:\|z\|<\alpha-1\}$.\\
We will prove that $G$ is analytic in $D_{\alpha+1}$, so we follow the same method as in (\ref{exp:22}) and (\ref{exp:23}) and obtain
\begin{align}
G(z)&=\sum_{n=0}^{\infty}\frac{z^{n}}{n!}\int_{0}^{\infty}\left(e^{i\lambda x}-1-i\lambda x\ind_{\{x<1\}}\right)x^{n}(e^{x}-1)^{-(\alpha+1)}dx\notag\\
&=\sum_{n=0}^{\infty}a_{n}z^{n},\notag
\end{align}
for $z\in D_{\alpha+1}$, where
\begin{equation}
a_{n}=\frac{1}{n!}\int_{0}^{\infty}\left(e^{i\lambda x}-1-i\lambda x\ind_{\{x<1\}}\right)x^{n}(e^{x}-1)^{-(\alpha+1)}dx\notag.
\end{equation}
In a similar way we prove that the function $I$ defined in (\ref{exp:12}) is an entire function.\\
We note that if $z\in D_{\alpha+1}$, then $\Re(-i\lambda+\alpha+1-z)>0$, and $\Re(\alpha+1-z)>0$, which implies that the function $P$ defined in (\ref{exp:14}) is also analytic in $D_{\alpha+1}$.\\
Finally since $G=i\lambda\tilde{a}+i\lambda I+P$ in $D_{\alpha-1}$, and both sides are analytic in $D_{\alpha+1}$, then $G=i\lambda\tilde{a}+i\lambda I+P$ in $D_{\alpha+1}$.\\
Therefore the first integral in the right-hand side of (\ref{exp:6}), $I_{1}$, is given by
\begin{equation}
I_{1}=G(\beta)=i\lambda \left(\tilde{a}_{\beta}+\frac{1}{\alpha-1}\right)+\Gamma(-\alpha)\left(\frac{\Gamma(-i\lambda+\alpha+1-\beta)}{\Gamma(-i\lambda+1-\beta)}-\frac{\Gamma(\alpha+1-\beta)}{\Gamma(1-\beta)}\right),\notag
\end{equation}
where
\begin{equation}\label{exp:30}
\tilde{a}_{\beta}=\tilde{a}+I(\beta)=\int_{0}^{1}\frac{xe^{-x}(1-e^{-(\alpha-\beta)x})}{(1-e^{-x})^{\alpha+1}}dx+\int_{0}^{1}e^{-x}\frac{1-x-e^{-x}}{(1-e^{-x})^{\alpha+1}}dx+\int_{1}^{\infty}\frac{e^{-x}}{(1-e^{-x})^{\alpha}}dx,
\end{equation}
for all $\beta<\alpha+1$. Now we compute the second integral in the right-hand side of (\ref{exp:6})
\begin{equation}
I_{2}=-i\lambda\tilde{b}_{\rho}+\int_{0}^{1}\frac{(u^{i\lambda}-1)u^{\rho_{1}-1}+i\lambda(u-1)}{(1-u)^{\alpha+1}}du,\notag
\end{equation}
where
\begin{equation}\label{exp:31}
\tilde{b}_{\rho}=\int_{0}^{1}\frac{xe^{-x}(1-e^{-(\alpha-\rho)x})}{(1-e^{-x})^{\alpha+1}}dx+\int_{0}^{1}e^{-x}\frac{1-x-e^{-x}}{(1-e^{-x})^{\alpha+1}}dx+\int_{1}^{\infty}\frac{e^{-x}}{(1-e^{-x})^{\alpha}}dx,
\end{equation}
and $\rho_{1}=\alpha+1-\rho$, hence following the same arguments used in the computation of $I_{1}$, we get
\begin{equation}
I_{2}=-i\lambda \left(\tilde{b}_{\rho}+\frac{1}{\alpha-1}\right)+\Gamma(-\alpha)\left(\frac{\Gamma(i\lambda+\alpha+1-\rho)}{\Gamma(i\lambda+1-\rho)}-\frac{\Gamma(\alpha+1-\rho)}{\Gamma(1-\rho)}\right),\notag
\end{equation}
for all $\rho<\alpha+1$. Therefore from the form of $I_{1}$ and $I_{2}$ we get
\begin{align}
\Psi_{L}(\lambda)&=-i\lambda\left(c_{+}\tilde{a}_{\beta}-c_{-}\tilde{b}_{\rho}+\frac{c_{+}-c_{-}}{\alpha-1}\right)-c_{+}\Gamma(-\alpha)\left((-i\lambda+1-\beta)_{\alpha}-(1-\beta)_{\alpha}\right)\notag\\
&-c_{-}\Gamma(-\alpha)\left((i\lambda+1-\rho)_{\alpha}-(1-\rho)_{\alpha}\right),\notag
\end{align}
for all $\beta, \rho<\alpha+1$.\\
$ii)$ Now we will compute the characteristic exponent when $\alpha=1$. In the following we assume that $\theta=0$ and that $c_{+}=c_{-}=1$. Since $\alpha=1$, the characteristic exponent of $X^{L}$ is given by
\begin{align}
\Psi_{L}(\lambda)&=-\bigg(c_{+}\int_{0}^{\infty}\left(e^{i\lambda x}-1-i\lambda x\ind_{\{x<1\}}\right)\frac{e^{\beta x}}{(e^{x}-1)^{2}}dx\notag\\
&+c_{-}\int_{-\infty}^{0}\left(e^{i\lambda x}-1-i\lambda x\ind_{\{x>-1\}}\right)\frac{e^{-\rho x}}{(e^{-x}-1)^{2}}dx\bigg).\notag
\end{align}
We will follow the same arguments used in the first part of the computation of the characteristic exponent in the case $\alpha\in(1,2)$. But to compute the two integrals in (\ref{exp:7}) we will need the following integral representation for the Digamma function (see \cite{GR})
\begin{equation}\label{exp:24}
\psi(z)=\int_{0}^{1}\frac{t^{z-1}-1}{z-1}dt-\mathcal{C},\qquad \textrm{for}\quad z\in\mathbb{C},
\end{equation}
where $\mathcal{C}$ is the Euler constant.\\
Now by making an integration by parts, using (\ref{exp:24}), and the recurrence relation for the Digamma function $\psi(z+1)=\psi(z)+z^{-1}$, we can express the first integral, for $\Re(z_{1})>1$, in the following form
\begin{align}\label{exp:17}
\int_{0}^{1}(u^{-i\lambda}-1)u^{z_{1}-2}(1-u)^{-1}du&=-\left(\int_{0}^{1}\frac{u^{-i\lambda+z_{1}-2}-1}{u-1}du-\int_{0}^{1}\frac{u^{z_{1}-2}-1}{u-1}du\right)\notag\\
&=\psi(z_{1}-1)-\psi(-i\lambda+z_{1}-1)\notag\\
&=\psi(z_{1})-\frac{1}{z_{1}-1}-\psi(-i\lambda+z_{1})+\frac{1}{-i\lambda+z_{1}-1},
\end{align}
As for the second integral in the right-hand side of (\ref{exp:7})
\begin{align}\label{exp:18}
\int_{0}^{1}(u^{-i\lambda+z_{1}-2}-1)(1-u)^{-1}du&=-\psi(-i\lambda+z_{1}-1)-\mathcal{C}\notag\\
&=\frac{1}{-i\lambda+z_{1}-1}-\psi(-i\lambda+z_{1})-\mathcal{C}.
\end{align}
So using (\ref{exp:17}) and (\ref{exp:18}), in (\ref{exp:7}) and recalling that $z_{1}=2-z$,
we get
\begin{align}\label{exp:19}
\int_{0}^{1}\frac{(u^{-i\lambda}-1)u^{z_{1}-1}+i\lambda(u-1)}{(1-u)^{2}}du=i\lambda(1-\mathcal{C})+(-i\lambda+1-z)\psi(-i\lambda+2-z)-(1-z)\psi(2-z)
\end{align}
if $\Re(z)<1$.\\
We note that (\ref{exp:19}) can be extended to the case where, $\Re(z)<2$, by the same arguments used in the case $\alpha\in(1,2)$, we only need to remark that the function $P:U\rightarrow\mathbb{C}$ defined by
\begin{equation}
P(z)=i\lambda(1-\mathcal{C})+(-i\lambda+1-z)\psi(-i\lambda+2-z)-(1-\beta)\psi(2-z)\notag,
\end{equation}
is analytic in the disc $D_{2}=\{z\in\mathbb{C}:\|z\|<2\}$. This implies that (\ref{exp:19}) is true for all $z\in D_{2}$, so in particular
\begin{equation}
I_{1}=G(\beta)=i\lambda(\tilde{a}_{\beta}+1-\mathcal{C})+(-i\lambda+1-\beta)\psi(-i\lambda+2-\beta)-(1-\beta)\psi(2-\beta)\notag,
\end{equation}
where
\begin{equation}
\tilde{a}_{\beta}=\int_{0}^{1}\frac{xe^{-x}(1-e^{-(1-\beta)x})}{(1-e^{-x})^{2}}dx+\int_{0}^{1}e^{-x}\frac{1-x-e^{-x}}{(1-e^{-x})^{2}}dx+\int_{1}^{\infty}\frac{e^{-x}}{(1-e^{-x})}dx,\notag
\end{equation}
for all $\beta<2$.\\
Finally by the same arguments used in the computation of $I_{1}$, we obtain
\begin{equation}
I_{2}=-i\lambda(\tilde{b}_{\rho}+1-\mathcal{C})+(i\lambda+1-\rho)\psi(i\lambda+2-\rho)-(1-\rho)\psi(2-\rho)\notag,
\end{equation}
where
\begin{equation}
\tilde{b}_{\rho}=\int_{0}^{1}\frac{xe^{-x}(1-e^{-(1-\rho)x})}{(1-e^{-x})^{2}}dx+\int_{0}^{1}e^{-x}\frac{1-x-e^{-x}}{(1-e^{-x})^{2}}dx+\int_{1}^{\infty}\frac{e^{-x}}{(1-e^{-x})}dx,\notag
\end{equation}
for all $\rho<2$.
Therefore from the form of $I_{1}$ and $I_{2}$ we get
\begin{align}
\Psi_{L}(\lambda)&=-i\lambda\left(c_{+}\tilde{a}_{\beta}-c_{-}\tilde{b}_{\rho}+(c_{+}-c_{-})(1-\mathcal{C})\right)\notag\\
&-c_{+}\bigg((-i\lambda+1-\beta)\psi(-i\lambda+2-\beta)-(1-\beta)\psi(2-\beta)\bigg)\notag\\
&-c_{-}\bigg((i\lambda+1-\rho)\psi(i\lambda+2-\rho)-(1-\rho)\psi(2-\rho)\bigg),\notag
\end{align}
for all $\beta, \rho<2$.\QED
Using the well known relationship between the Laplace and the characteristic exponents, we obtain:
\begin{corollary} Let $X^{L}$ be a Lamperti stable process with characteristics $(\alpha,f,\sigma, \theta)$.
\begin{itemize}
\item[i)] Let $\alpha\in(0,1)$ and suppose that $X^L$ is a Lamperti stable subordinator, then its Laplace exponent is given by
 \[
\Phi_{L}(\lambda)=\mathbf{d}\lambda-c_+\Gamma(-\alpha)\bigg((\lambda+1-\beta)_{\alpha}-(1-\beta)_{\alpha}\bigg), \qquad \lambda\geq 0,
\]
where $\mathbf{d}\geq0$.
\item[ii)] Let $\alpha\in(1,2)$ and suppose that $X^{L}$ has no positive jumps then its Laplace exponent is given by
\begin{equation}\label{snlexp}
\Phi_{L}(z)=-\tilde{\theta}\lambda+c_-{\Gamma(-\alpha)}\bigg((\lambda+1-\rho)_{\alpha}-(1-\rho)_{\alpha}\bigg), \qquad \lambda\geq0,\notag
\end{equation}
where $\tilde{\theta}$ is given by (\ref{exp:40}).
\end{itemize}
\end{corollary}

\begin{remark}
This Corollary has, as particular cases, the two recent results found in:
\begin{itemize}
\item[i)] Corollary 2, and Lemma 4, in \cite{ckp}, where the result is obtained by means of the Lamperti transformation.
\item[ii)] Proposition 3.1 in \cite{pp}, where the Laplace exponent is obtained using special functions. The ideas in \cite{pp} as well as in \cite{cc}, inspired parts of the proof if Theorem 1, specifically the decomposition (\ref{exp:11}).
\end{itemize}
\end{remark}

Now we turn our attention to other group of properties. Let $H=(H_t,t\geq 0)$ be the increasing ladder height process of
$X^{L}$ (see chapter VI in \cite{be}) and
$\widehat{H}=(\widehat{H}_t, t\geq 0)$, its decreasing ladder height
process. Denote by $k$ and $\hat{k}$ for the characteristic
exponents of  $H$ and $\widehat{H}$, which are subordinators, and
suppose that $X^{L}$ drifts to $-\infty$ and $\nu_{\sigma}^{\alpha,
f}(0,\infty)>0.$ Under this hypothesis, the process $H$ is a killed
subordinator and we denote by $\Pi_H$ for its L\'evy measure.  The
following result give us a relation between $\nu^{\alpha,
f}_{\sigma}$ and $\Pi_H$.
\begin{proposition} Let $X^{L}$ be a Lamperti stable process with positive jumps and characteristics $(\alpha, f)$ such that it drifts to $-\infty$. Then,  the tail of the L\'evy measure of $H$, its increasing ladder height process, belongs to $\mathcal{L}^{(\alpha+1-\beta)}$ and
\[
\nu^{\alpha,f}_\sigma(u,\infty)\sim
\hat{k}(-i(\alpha+1-\beta))\Pi(u,\infty)\qquad \textrm{as }\quad
u\to\infty.
\]
\end{proposition}
\noindent {\it Proof:} The proof follows directly from Proposition
5.3 in \cite{kkm} and proposition 4.\QED

We finish this section with some properties of  Lamperti stable
processes with no positive jumps.
\begin{proposition}
Let $X^{L}$ be a Lamperti stable process with no positive jumps and
characteristics $(\alpha, \rho,\sigma,\theta)$, such that $\tilde{\theta}=0$ in (\ref{exp:40}). Then,
\begin{itemize}
\item[i)] there exist $\rho_0\in (1,2)$ such that $X_{LS}$ drifts to $\infty$, oscillates or drifts to $-\infty$ according as $\rho\in(-\infty,\rho_0)$, $\rho=\rho_0$ or $\rho\in(\rho_0,\alpha+1)$.
\item[ii)] for $\rho\in(\rho_0,\alpha+1)$, we have that there exist $\lambda>0$ such that
\begin{equation}\label{estsup}
\mathbf{P}_0\Big(S^{LS}_\infty>x\Big)\sim \frac{c}{\lambda
k}e^{-\lambda x}, \quad \textrm{as}\quad x\to \infty,
\end{equation} where $S^{LS}_{\infty}=\sup_{t\geq 0} X^{L}_t$,
$c=-\log \mathbf{P}_0(H_1<\infty)$, $k=\mathbf{E}_0(H_1e^{\lambda
H_1};H_1<\infty)$ and $H$ is the increasing ladder height process.
\item[iii)] for $\rho\in(\rho_0,\alpha+1)$, we have that there exist $\lambda>0$ such that
\begin{equation}\label{estint}
\mathbf{P}_0\Big(I(X^{L})>x\Big)\sim Kx^{-\lambda}, \quad
\textrm{as}\quad x\to \infty,
\end{equation}
where $K$ is a positive constant  and
\[
I(X^{L})=\int_0^{\infty}\exp\Big\{X^{L}_t\Big\}\ud t.
\]
\item[iv)] the process $X^{L}$ has increase times.
\item[v)] the process $X^{L}$ satisfies the Spitzer's condition at $\infty$, i.e.
\[
\lim_{t\to \infty}\frac{1}{t}\int_{0}^{t}\mathbf{P}(X^{L}_s\geq
0)\ud s =1/\alpha  \quad \textrm{as}\quad x\to \infty.
\]
\item[vi)] the process $X^{L}$ satisfies the following law of the iterated logarithm
\begin{equation}\label{lil}
\limsup_{x\to 0}\frac{X^{L}_t \Phi^{-1}_{L}(t^{-1}\log|\log
t|)}{\log|\log t|}=\alpha^{-\alpha}(1-\alpha)^{\alpha-1} \qquad a.s.,
\end{equation}
where $\Phi^{-1}_{L}$ denotes the right-continuous inverse of
$\Phi^{-1}_{L}$.
\end{itemize}
\end{proposition}
\noindent {\it Proof:} $(i)$ We know that in this case $\alpha\in (1,2)$, so from Corollary VII.2 in
 \cite{be}, the process $X^{L}$ drifts to $+\infty$, oscillates or
 drifts to $-\infty$ according as $\Phi_{L}'(0^+)$ is positive, zero
 or negative.  Hence, from the Laplace exponent of $X^L$ we have, using the recursion formula
 for the Gamma and Digamma functions, the following
\begin{align}\label{d0L}
\Phi_{L}'(0^+)&=c_-\Gamma(-\alpha)(1+\alpha-\rho)_{\alpha}(\psi(1-\rho+\alpha)-\psi(1-\rho)),\notag\\
&=c_-\Gamma(-\alpha)\frac{\Gamma(1+\alpha-\rho)}{\Gamma(3-\rho)}\bigg((2-\rho)(1-\rho)((\psi(1-\rho+\alpha)-\psi(1-\rho))+3-2\rho\bigg),\notag\\
&:=g(\rho).
 \end{align}
 We have from (\ref{d0L}) that $g(1)<0$, and $g(2)>0$. On the other hand, in the interval $(1,2)$, the function $g$ is continuous and decreasing which implies that there exist $\rho_0\in (1,2)$ such that $g(\rho_0)=0$. Thus, we  deduce that $X^{L}$ drifts to $\infty$, oscillates or drifts to $-\infty$ according as $\rho\in(-\infty,\rho_0)$, $\rho=\rho_0$ or $\rho\in(\rho_0,\alpha+1)$.\\
 $(ii)$ First, note that any L\'evy process with no positive jumps which drifts to $-\infty$ satisfies that its Laplace exponent has a strictly positive root. Hence for a Lamperti stable process with no positive jumps and with  $\rho\in (\rho_0,\alpha+1)$, there exists $\lambda>0$
 such that
 \[
 \mathbf{E}_0\Big(\exp\{\lambda X^{L}_1\}\Big)=1,
 \]
 i.e. that $X^L$ satisfies the Cram\'er condition. Thus,   the main result in \cite{bd} give us the sharp estimate in (\ref{estsup}).\\
 $(iii)$ First note that $X^{L}$ is not arithmetic and that under our
 assumptions  the Cram\'er condition is satisfied for some
 $\lambda>0$. Hence from  Lemma 4 in \cite{ri}, we get the sharp
 estimate (\ref{estint})
 for the exponential functional $I(X^{L})$.\\
$(iv)$ Here, we need the following estimate of the Pochhammer symbol
(see for instance \cite{Le}),
\begin{equation}\label{estle}
(\lambda+1-\rho)_{\alpha}\sim \lambda^{\alpha}\qquad
\textrm{as}\quad \lambda \to \infty.
\end{equation}
From Corollary VII.9 and Proposition VII.10 in \cite{be} we know
that $X^{L}$ has increase times if
\[
\int^{\infty} \lambda^{-3}\Phi_{L}(\lambda)\ud \lambda<\infty,
\]
which in our case is satisfied since from (\ref{estle}), we have
\begin{equation}\label{estle1}
 \Phi_{L}(\lambda)\sim c_-\Gamma(-\alpha)\lambda^{\alpha}\qquad \textrm{as}\quad \lambda \to \infty.
\end{equation}
$(v)$ From (\ref{estle1}), we  see that $\Phi_{L}$ is regularly varying at $\infty$ with index $\alpha$. Hence, the statement follows from Proposition VII.6 in \cite{be}. \\
$(vi)$ Since $\Phi_{L}$ is regularly varying at $\infty$ with index
$\alpha$, we have that its right-continuous inverse $\Phi_{L}^{-1}$
is regularly varying at $\infty$ with index $1/\alpha$ which
corresponds to the Laplace exponent of the first passage time of
$X^{L}$ (which is a subordinator). Therefore, from Theorem III.11 in
\cite{be} we deduce the law of the iterated logarithm
(\ref{lil}).\QED
\section{Short and long time behaviour.}
Motivated by  the works of Rosi\'nski \cite{ro} and Houdr\'{e} and
Kawai \cite{hok}, we study
 the short and long time behavior of  Lamperti stable
processes. In particular, we will show that this class of  processes
share with the tempered and layered stable processes, the
peculiarity that in short
time they behave like  stable processes.

 The convergence in distribution of processes, considered in this section, is in the functional sense, i.e.  in the sense of the weak convergence of the laws of the processes on the Skorokhod space and will be denoted by $``\overset{d}{\rightarrow}"$.
\begin{proposition}Let $X^{L}$ be a Lamperti stable process with characteristics $(\alpha, f,\sigma,0)$ and
\begin{equation}
\eta_{\alpha}= \left\{
 \begin{array}{ll}
  0&\hbox{if $\alpha=1$,}\notag\\
 \displaystyle\int_{S^{d-1}}\xi\sigma(\ud\xi)\displaystyle\int_{0}^{1}re^{f(\xi) r}(e^{r}-1)^{-(\alpha+1)}\ud r&\hbox{if $\alpha\in(0,1)$},\notag\\
\displaystyle\int_{S^{d-1}}\xi\sigma(\ud\xi)\displaystyle\int_{1}^{\infty}re^{f(\xi) r}(e^{r}-1)^{-(\alpha+1)}\ud r&\hbox{if $\alpha\in(1,2)$.}\notag\\
 \end{array}
\right.
\end{equation}
Then,
\begin{equation}
\Big(h^{-1/\alpha}\big(X^{L}_{ht}-ht\eta_{\alpha}\big),t>0\Big)\overset{d}{\rightarrow}(X_{t},t>0)\quad
as\quad h\rightarrow 0,\notag
\end{equation}
where $(X_{t},t>0)$ is a stable process of index $\alpha$.
\end{proposition}
\noindent{\it Proof:} The proof is similar to that of the short time
behaviour of layered stable process, since for each $\xi\in S^{d-1}$
\[
e^{f(\xi) r}(e^{r}-1)^{-(\alpha+1)}\sim r^{-(\alpha+1)}\quad as\quad
r\rightarrow0.
\]
Thus, we follow the proof of Theorem 3.1 in \cite{hok} with
\[
q(\xi,r)=e^{f(\xi) r}(e^{r}-1)^{-(\alpha+1)}
\]
and the desired result is obtained.\QED
\begin{theorem} Let $X^{L}_{t}$ be a Lamperti stable process with characteristics $(\alpha, f,\sigma,0)$ and
\[
\eta_{\alpha}=-\int_{S^{d-1}}\xi\sigma(\ud\xi)\int_{1}^{\infty}re^{f(\xi)
r}(e^{r}-1)^{-(\alpha+1)}\ud r.
\]
Then,
\begin{equation}\label{ltb}
\Big(h^{-1/2}\big(X^{L}_{ht}-ht\eta_{\alpha}\big),t>0\Big)\overset{d}{\rightarrow}(W_{t},t>0)\quad
as\quad h\rightarrow \infty,
\end{equation}
where $(W_t,t>0)$ is a centered Brownian motion with covariance
matrix
\[
\int_{\scriptsize \R^{d}_{0}}xx'\nu^{\alpha,f}_{\sigma}(\ud x).
\]
\end{theorem}
\noindent{\it Proof:} According to  a standard result on the
convergence of processes with independent increments due to
Skorokhod (see for instance Theorem  15.17 of Kallenberg
\cite{kal}), the functional convergence (\ref{ltb}) holds if and
only if
\[
h^{-1/2}\big(X^{L}_{h}-h\eta_{\alpha}\big)\overset{d}{\rightarrow}W_{1}\quad
as\quad h\rightarrow \infty.
\]
Now, we introduce the following transform for positive measures, for
any $r>0$
\[
(T_r\nu)(B)=\nu(r^{-1}B)\qquad\textrm{for}\quad B\in
\mathcal{B}(\R^d).
\]
Note that the random variable $h^{-1/2}X_{h}^{LS}$ is infinitely
divisible  and since it has finite first moment, we may rewrite its
characteristic exponent as follows;
\[
ih\int_{\scriptsize\R^{d}_0}\langle y,x\rangle\ind_{\{\|x\|\geq
1\}}(T_{h^{-1/2}}\nu^{\alpha,f}_{\sigma})(\ud
x)-h\int_{\scriptsize\R^{d}_{0}}\Big(e^{i\langle
y,x\rangle}-1-i\langle
y,x\rangle\ind_{\{\|x\|\leq1\}}\Big)(T_{h^{-1/2}}\nu^{\alpha,f}_{\sigma})(\ud x).
\]
Hence, from Theorem 15.14 of Kallenberg \cite{kal} we only need  to
check the following convergences as $h$ increases:
\begin{itemize}
\item[a)] $h(T_{h^{-1/2}}\nu^{\alpha,f}_{\sigma})$ converges vaguely towards $0$ on $\R^d_0$,
\item[b)] for each $k>0$, \hspace{1cm}$h\displaystyle\int_{\|x\|\leq\kappa}xx'(T_{h^{-1/2}}\nu^{\alpha,f}_{\sigma})(\ud x)\rightarrow\displaystyle\int_{\mathbb{R}^{d}_{0}}xx'\nu^{\alpha,f}_{\sigma}(\ud x),$
\item[c)] for each $k>0$, \hspace{1cm}$h\displaystyle\int_{\|x\|\geq k}x(T_{h^{-1/2}}\nu^{\alpha,f}_{\sigma})(\ud x)\rightarrow0.$
\end{itemize}
We first prove $(a)$ or equivalently
\begin{equation}\label{a:11}
\lim_{h\to\infty}
\int_{\scriptsize\R^{d}_{0}}g(x)h(T_{h^{-1/2}}\nu^{\alpha,f}_{\sigma})(\ud
x)=0
\end{equation}
for all bounded continuous functions $g:\R_{0}^{d}\rightarrow\R$
vanishing in a neighborhood of the origin. Let $g$ be such a
function satisfying that $|g|\leq C$, and that for some $\delta>0$,
$g(x)\equiv0$ on $\{x\in \R^{d}_{0}:\|x\|<\delta\}$. Let $\gamma:=\sup_{\xi\in S^{d-1}}f(\xi)$, then we have
\begin{align}\label{a:10}
\bigg|h\int_{\scriptsize\R^{d}_{0}}g(x)\,&(T_{h^{-1/2}}\nu^{\alpha,f}_{\sigma})(\ud x)\bigg|\notag\\
&\leq
h^{1+1/2}\int_{S^{d-1}}\sigma(\ud\xi)\int_{0}^{\infty}|g(r\xi)|e^{rf(\xi)
h^{1/2}}(e^{rh^{1/2}}-1)^{-(\alpha+1)}\ud r\notag\\
&=\int_{S^{d-1}}\sigma(\ud\xi)\int_{\delta}^{\infty}|g(r\xi)|\frac{(rh^{1/2})^{3}}{r^{3}}e^{rh^{1/2}\gamma}(e^{rh^{1/2}}-1)^{-(\alpha+1)}\ud r.\notag\\
\end{align}
On the other hand, since $\gamma<\alpha+1$ it follows
\begin{equation}
\lim_{r\to\infty}r^{3}\frac{e^{r\gamma}}{(e^{r}-1)^{\alpha+1}}=0,\notag
\end{equation}
then for  $\epsilon>0$ sufficiently small, there exist  $M>0$ such
that for all $r\geq M$
\begin{equation}
r^{3}e^{rf(\xi)}(e^{r}-1)^{-(\alpha+1)}<\epsilon.\notag
\end{equation}
Since $r>\delta$, we may take $h>\left(\frac{M}{\delta}\right)^{2}$
in (\ref{a:10}) and obtain
\begin{align}
\bigg|h\int_{\scriptsize\R^{d}_{0}}g(x)\,(T_{h^{-1/2}}\nu^{\alpha,f}_{\sigma})(\ud x)\bigg|&<\epsilon\int_{S^{d-1}}\sigma(d\xi)\int_{\delta}^{\infty}|g(r\xi)|\frac{1}{r^{3}}\ud r\notag\\
&\le\epsilon
C\int_{S^{d-1}}\sigma(d\xi)\int_{\delta}^{\infty}\frac{1}{r^{3}}\ud
r.\notag
\end{align}
Note that the last integral in the right-hand side of the above
inequality is finite and therefore the convergence (\ref{a:11})
follows.\\
Next, we prove part $(b)$. First note that
$\displaystyle\int_{\scriptsize\R^{d}_{0}}\|x\|^2\nu^{\alpha,f}_{\sigma}(\ud
x)$ is finite. This follows by similar arguments as those used in
proposition 1. This implies that the integral
$\displaystyle\int_{\scriptsize\R^{d}_{0}}xx'\nu^{\alpha,f}_{\sigma}(\ud
x)$ is well defined. Now take $k>0$ fixed, and note that
\[
h\int_{\{\|x\|\leq k\}}xx'(T_{h^{-1/2}}\nu^{\alpha,f}_{\sigma})(\ud
x)=\int_{\{\|x\|\leq h^{1/2}k\}}xx'\nu^{\alpha,f}_{\sigma}(\ud
x)\to\int_{\scriptsize\R^{d}_{0}}xx'\nu^{\alpha,f}_{\sigma}(\ud x),
\]
as $h$ goes to $\infty$, which proves part $(b)$.\\
Finally, we  consider $k>0$ and recall that $\gamma=\sup_{\xi\in
S^{d-1}}f(\xi)$, then
\begin{align}
\bigg\|h\int_{\{\|x\|\geq k\}}&z(T_{h^{-1/2}}\nu^{\alpha,f}_{\sigma})(\ud z)\bigg\|\notag\\
&=\left\|h^{1+1/2}\int_{S^{d-1}}\xi\sigma(\ud\xi)\int_{k}^{\infty}re^{rf(\xi)
h^{1/2}}(e^{rh^{1/2}}-1)^{-(\alpha+1)}\ud r\right\|\notag\\
&\leq(1-e^{-kh^{1/2}})^{-(\alpha+1)}\left\|h^{1+1/2}\int_{S^{d-1}}\xi\sigma(\ud\xi)\int_{k}^{\infty}re^{r
h^{1/2}(\gamma-(\alpha+1))}\ud r\right\|\notag\\
&=\frac{e^{-kh^{1/2}(\alpha+1-\gamma)}}{(1-e^{-kh^{1/2}})^{\alpha+1}}\left(\frac{hk}{\alpha+1-\gamma}-\frac{h^{1/2}}{(\alpha+1-\gamma)^2}\right)\left\|\int_{S^{d-1}}\xi\sigma(d\xi)\right\|,\notag
\end{align}
which goes to $0$ as $h\rightarrow\infty$ since $\gamma<\alpha+1$.
This completes the proof.\QED Let us apply the above results to the
special cases treated in the introduction. In particular when we
start with a stable process $(X,\mathbf{P}_x),x>0$, of index
$\alpha$,  applying the result in short time behavior  after various
transformations we return to this initial process. Recall that
associated to the stable process  three L\'evy processes are
obtained via the Lamperti representation of pssMp:
$\xi^*,\xi^\uparrow,\xi^\downarrow$. Then the normalization of any
of them according to proposition 10, converges weakly in the space
of Skorokhod to the original stable process $X$, i.e.
\begin{align*}
X &\overset{kill}\longrightarrow X^* \overset{LT}\longrightarrow X^{L} \overset{norm}\longrightarrow X^{L}_h  \overset{d}{\rightarrow} X \qquad \textrm{as}\quad h\rightarrow 0 \\
X &\overset{kill}\longrightarrow X^* \overset{DT}\longrightarrow
X^C\overset{LT}\longrightarrow X^{L} \overset{norm}\longrightarrow
X^{L}_h  \overset{d}{\rightarrow} X \quad \textrm{as}\quad
h\rightarrow 0
\end{align*}
where $kill$, $LT$, $DT$ and $norm$ means killing , the Lamperti
representation of pssMp, Doob-transform or conditioning, and
normalization of a given process, respectively. Moreover $X^C$ is
the conditioned process (to be positive or to hit $0$ continuously),
$X^{L}$ stands for any of the Lamperti stable processes
$\xi^\uparrow,\xi^\downarrow$ and $\xi^*$, and $X^{L}_h$ is  the
normalization of each of them given in proposition 10. In the same
spirit we could also write, using theorem 2,
\begin{align*}
X &\overset{kill}\longrightarrow X^* \overset{LT}\longrightarrow X^{L} \overset{norm}\longrightarrow X^{L}_h  \overset{d}{\rightarrow} W \qquad \textrm{as}\quad h\rightarrow \infty, \\
X &\overset{kill}\longrightarrow X^* \overset{DT}\longrightarrow
X^C\overset{LT}\longrightarrow X^{L} \overset{norm}\longrightarrow
X^{L}_h  \overset{d}{\rightarrow} W \quad \textrm{as}\quad
h\rightarrow \infty,
\end{align*}
where $W$ is a centered brownian motion. \\

The final result of this section follows the line of reasoning  of
last remark but uses additional tools that we shall briefly
introduce. In \cite{clu} the convergence in the Skorokhod space is
studied in relation to the second Lamperti transformation $(LT_2)$,
i.e. the one that transforms  L\'evy processes with no negative
jumps to continuous state branching processes. The problem of
explosions is difficult to handle in this metric so the authors
consider another metric $d_{\infty}$ on the Skorohod space which
given by
$$d_{\infty}(f,g)=1\land\inf_{\lambda\in\Lambda_\infty}\|f-g\circ \lambda\|_\infty\vee \|\lambda-I\|_\infty.$$
and where $\Lambda_{\infty}$ is the set of increasing
homeomorphisms of  $[0,\infty)$ into itself.
 According to the authors  the convergence in this metric implies it in the usual Skorohod metric. Two of the main results in \cite{clu}  say that  the Lamperti transform $LT_2$,  is continuous with this new metric (proposition 4) and that a sequence $(Y^{*,n})$ of stopped L\'evy processes with no negative jumps converges in this new metric  towards $Y$, a stopped L\'evy process when the sequence of the associated Laplace exponents of  $(Y^{*,n})$ converges towards the associated Laplace exponent of $Y$ (proposition 5).  Therefore, a combination of the results mentioned above and proposition 8 give us the following corollary.

\begin{corollary}
Let $X$ and $X^{L}$ be a stable proces of index $\alpha$ with no
negative jumps and a Lamperti stable processes with no negative
jumps with characteristics $(\alpha,\beta)$ which does not drift
towards $+\infty$, respectively. Let $Y_h = LT_2(X^{L}_h )$ and $Y =
LT_2(X)$. Then
$$Y_h  \overset{d}{\rightarrow} Y\qquad \textrm{as} \quad h\to 0.$$
\end{corollary}

\section{Absolute continuity with respect to stable processes}
We  showed that in short time  a Lamperti-stable process behaves
like a stable process, now following Rosi\'nski \cite{ro} we will
relate the law of both processes. In other words, we will find a
probability measure under which the law of  a Lamperti stable
process with characteristics $(\alpha,f,\sigma)$ is the same that
the law of the short time limiting stable process with index
$\alpha$.
\begin{theorem}
Let $P$ and $Q$ be two probability measures on
$(\Omega,\mathcal{F})$ and such that under $P$ the canonical process
$(X_{t},t\geq0)$ is a Lamperti stable process with characteristics
$(\alpha,f,\sigma,a)$, while under $Q$ it is a stable  process with index
$\alpha$ with linear term $b$. Let $(\mathcal{F}_t)$ be the
canonical filtration, and assume that $f\in L^{2}(S^{d-1},\mathbb{B}(S^{d-1}),\sigma)$.
Then
\begin{itemize}
\item[i)] $P|_{\mathcal{F}_{t}}$ and $Q|_{\mathcal{F}_{t}}$ are
mutually absolutely continuous for every $t>0$ if and only if
\begin{equation}
a-b=\notag
 \left\{
              \begin{array}{ll}
               \displaystyle\int_{S^{d-1}}\xi\sigma(\ud\xi) \displaystyle\int_{0}^{1}re^{rf(\xi)}(e^{r}-1)^{-(\alpha+1)}\ud r, & \hbox{if $\alpha\in (0,1)$,}\notag\\
                 \displaystyle\int_{S^{d-1}}\xi\sigma(\ud\xi) \displaystyle\int_{0}^{1}r(e^{rf(\xi)}(e^{r}-1)^{-(\alpha+1)}-r^{-(\alpha+1)})\ud r,& \hbox{if
$\alpha=1$,}\notag\\
                 \displaystyle\int_{S^{d-1}}\xi\sigma(\ud\xi) \displaystyle\int_{0}^{1}r(e^{rf(\xi)}(e^{r}-1)^{-(\alpha+1)}-r^{-(\alpha+1)})\ud r\notag\\
                -\displaystyle\int_{S^{d-1}}\xi\sigma(\ud\xi) \displaystyle\int_{1}^{\infty}r^{-(\alpha+1)}\ud r,& \hbox{if $\alpha\in(1,2)$.}
              \end{array}
\right.
\end{equation}
\item[ii)] For each $t>0$,
\begin{equation}
\frac{\ud Q}{\ud P}\bigg|_{\mathcal{F}_{t}}=e^{U_{t}},\notag
\end{equation}
where $(U_{t},t\geq0)$ is a L\'{e}vy process defined on
$(\Omega,\mathcal{F},P)$ by
\begin{align}
U_{t}=\lim_{\epsilon\downarrow0}\underset{\{s\in(0,t]:\|\Delta
X_{s}\|>\epsilon\}}{\sum}&\bigg[\Big(e^{\|\Delta X_{s}\|f(\Delta
X_{s})}(e^{\|\Delta X_{s}\|}-1)^{-(\alpha+1)}\|\Delta
X_{s}\|^{\alpha+1}\Big)\notag\\
&-t(\nu^{\alpha,
f}_{\sigma}-\Pi)\Big(\{z\in\mathbb{R}^{d}_{0}:\|z\|>\epsilon\}\Big)\bigg].\notag
\end{align}
In the above right hand side, the convergence holds $P$-a.s.
uniformly in $t$ on every interval of positive length.
\end{itemize}
\end{theorem}
\noindent{\it Proof:} From  Theorem 33.2 in Sato \cite{sa}, we only
need to verify that
\begin{equation}
\int_{\scriptsize \R_0^{d}}(e^{\varphi(x)/2}-1)^{2}\Pi(\ud x)<\infty,\notag
\end{equation}
where $\varphi:\R_0^{d}\rightarrow\R$ is defined by
\begin{equation}
\frac{\ud\nu^{\alpha,f}_{\sigma}}{\ud \Pi} (x)=e^{\varphi(x)}.\notag
\end{equation}
In particular, we have
$\varphi(r\xi)=\log\left(e^{rf(\xi)}(e^{r}-1)^{-(\alpha+1)}r^{\alpha+1}\right)$.
Thus, we need   to check
\begin{equation}\label{a:12}
\int_{S^{d-1}}\sigma(\ud
\xi)\int_{0}^{\infty}\left[\left(\frac{e^{rf(\xi)}r^{(1+\alpha)}}{(e^{r}-1)^{(\alpha+1)}}\right)^{1/2}-1\right]^{2}\frac{1}{r^{1+\alpha}}\ud
r<\infty
\end{equation}
By Taylor expansion and the Lagrange form for the residual, we have
$(e^{r}-1)=re^{r\theta_{r}},$ where $\theta_{r}\in(0,1)$. This
implies
\begin{equation}\label{a:15}
\frac{e^{rf(\xi)}r^{(1+\alpha)}}{(e^{r}-1)^{(\alpha+1)}}=e^{r(f(\xi)-\theta_{r}(\alpha+1))}.
\end{equation}
Now, noting that $f(\xi)-(\alpha+1)\leq
f(\xi)-\theta_{r}(\alpha+1)\leq f(\xi)$, it follows
\begin{equation}
e^{r(f(\xi)-(\alpha+1))/2}-1\leq
e^{r(f(\xi)-\theta_{r}(\alpha+1))/2}-1\leq e^{rf(\xi)/2}-1,\notag
\end{equation}
and since $f(\xi)\leq\gamma=\sup_{\xi\in S^{d-1}}f(\xi)$, we have
\begin{align}\label{a:13}
\Big(e^{r(f(\xi)-\theta_{r}(\alpha+1))/2}-1\Big)^{2}&\leq\Big(e^{r(f(\xi)-(\alpha+1))/2}-1\Big)^{2}\vee\Big(e^{rf(\xi)/2}-1\Big)^{2}
\end{align}
Using a Taylor expansion again and (\ref{a:13}), it is clear that
there exists a constant $R>0$ such that if $r<R$, then
\begin{equation}\label{a:14}
\Big(e^{r(f(\xi)-\theta_{r}(\alpha+1))/2}-1\Big)^{2}\leq K(f^{2}(\xi)+1)r^{2},
\end{equation}
where $K$ is a positive constant. Hence from (\ref{a:15}) and
(\ref{a:14}), it follows that
\begin{align}
\int_{S^{d-1}}\sigma(\ud\xi)\int_{0}^{R}&\left[\left(\frac{e^{rf(\xi)}r^{(1+\alpha)}}{(e^{r}-1)^{(\alpha+1)}}\right)^{1/2}-1\right]^{2}\frac{1}{r^{1+\alpha}}\ud r\notag\\
&\leq K\left(\sigma(S^{d-1})+\int_{S^{d-1}}f^2(\xi)d\xi\right)\int_{0}^{R}\frac{r^{2}}{r^{1+\alpha}}\ud
r,\notag
\end{align}
which is finite because $\alpha\in(0,2)$ and $f\in L^{2}(S^{d-1},\mathbb{B}(S^{d-1}),\sigma)$. In the case when $r>R$,  we have
\begin{align}
\int_{S^{d-1}}&\sigma(\ud\xi)\int_{R}^{\infty}\left[\left(\frac{e^{rf(\xi)}r^{(1+\alpha)}}{(e^{r}-1)^{(\alpha+1)}}\right)^{1/2}-1\right]^{2}\frac{1}{r^{1+\alpha}}\ud r\notag\\
&\leq4\left((1-e^{-R})^{-(\alpha+1)}\int_{S^{d-1}}\sigma(\ud\xi)\int_{R}^{\infty}e^{r(f(\xi)-(\alpha+1))}\ud r+\sigma(S^{d-1})\int_{R}^{\infty}\frac{1}{r^{1+\alpha}}\ud r\right)\notag\\
&\leq4\sigma(S^{d-1})\left((1-e^{-R})^{-(\alpha+1)}\int_{R}^{\infty}e^{r(\gamma-(\alpha+1))}\ud
r+\int_{R}^{\infty}\frac{1}{r^{1+\alpha}}\ud r\right),\notag
\end{align}
which is also finite because $\gamma<\alpha+1$. Therefore (\ref{a:12}) follows.

The proof of the second statement of the Theorem follows directly
from Theorem 33.2 of Sato \cite{sa}.\QED Here, we follow the same
notation as in Theorem 3. Note that under the conditions of Theorem
4.1 in \cite{hok}, if $R$ is another probability measure on
$(\Omega, \mathcal{F})$ under which the canonical process $X=(X_t,
t\geq 0)$ is a layered stable process, we have that
$R|_{\mathcal{F}_{t}}$ and $Q|_{\mathcal{F}_{t}}$ are mutually
absolutely continuous for every $t>0$. From our previous result, we
obtain the corresponding result for Lamperti stable processes, i.e.
that  $R|_{\mathcal{F}_{t}}$ and $P|_{\mathcal{F}_{t}}$ are mutually
absolutely continuous for every $t>0$. Similar result holds for the
tempered stable processes, see Theorem 4.1 in \cite{ro}.
\section{Series representations of Lamperti stable process}
In this section, we establish a series representation for Lamperti
stable processes which allow us to generate some of their sample
paths. To this end, we will
use the LePage's method found in \cite{lp}. We first introduce the following sequences of mutually independent random
variables defined in $[0,T]$. Let $\{\Gamma_{i}\}_{i\geq1}$ be a
sequence of of partial sums of iid standard exponential random
variables, $\{U_{i}\}_{i\geq1}$  be a sequence of uniform random
variables on $[0,T]$, and  let $\{V_{i}\}_{i\geq1}$ be a sequence of
iid random variables in $S^{d-1}$ with common distribution
$\sigma(d\xi)/\sigma(S^{d-1})$. In order to use the LePage's method,
we consider the following function
$\rho^{-1}:(0,\infty)\times S^{d-1}\to \R_+$ given by
\begin{equation}
\rho^{-1}(u,\xi):=\inf\Big\{x>0:\rho([x,\infty),\xi)<u\Big\},\notag
\end{equation}
where
\begin{equation}
\rho([x,\infty),\xi)=\int_{x}^{\infty}e^{f(\xi)r}(e^{r}-1)^{-(\alpha+1)}\ud
r.\notag
\end{equation}
Now, let $\{c_{i}\}_{i\geq1}$ be a sequence of constants defined as
follows,
\begin{equation}
c_{i}=\int_{i-1}^{i}\mathbf{E}\Big(\rho^{-1}(s/T,V_{1})V_{1}\ind_{\{\rho^{-1}(s/T,V_{1})\leq1\}}\Big)ds.\notag
\end{equation}
Then from Theorem 5.1 in \cite{ro2}, the process
\begin{equation}
\Bigg(\sum_{i=1}^{\infty}\Big(\rho^{-1}(\Gamma_{i}/T,V_{i})V_{i}\ind_{\{U_{i}\leq
t\}}-c_{i}\frac{t}{T}\Big),t\in[0,T]\Bigg),\notag
\end{equation}
converges uniformly a.s. towards a Lamperti stable process
with characteristics $(\alpha,f,\sigma)$ and linear term $\theta=0$ (in the L\'evy-Khintchine formula). In particular when $f(\xi)=1$, we have  that
\begin{equation}
\rho^{-1}(u,\xi)=\ln(1+(\alpha u)^{-1/\alpha}),\notag\\
\end{equation}
hence the series representation for a Lamperti stable L\'{e}vy
process $X^{L}$ with characteristics $(\alpha,1)$, is as follows
\begin{equation}
X^{L}_{t}\overset{d}{=}\sum_{i=1}^{\infty}\Bigg(\ln\bigg(1+\bigg(\frac{\alpha\Gamma_i}{T}\bigg)^{-1/\alpha}\bigg)V_{i}\ind_{\{U_{i}\leq
t\}}-c_{i}\frac{t}{T}\Bigg)\notag
\end{equation}
where
\begin{equation}
c_{i}=\mathbf{E}\Big(V_{1}\Big)\int_{i-1}^{i}\ln\bigg(1+\bigg(\frac{\alpha
s}{T}\bigg)^{-1/\alpha}\bigg)\ind_{\{\ln(1+(\alpha
sT^{-1})^{-1/\alpha})\leq1\}}\ud s.\notag
\end{equation}
Let us observe below some sample paths of this particular Lamperti stable process  generated via the series representation.
\begin{figure}
\centering
\includegraphics[scale=.37]{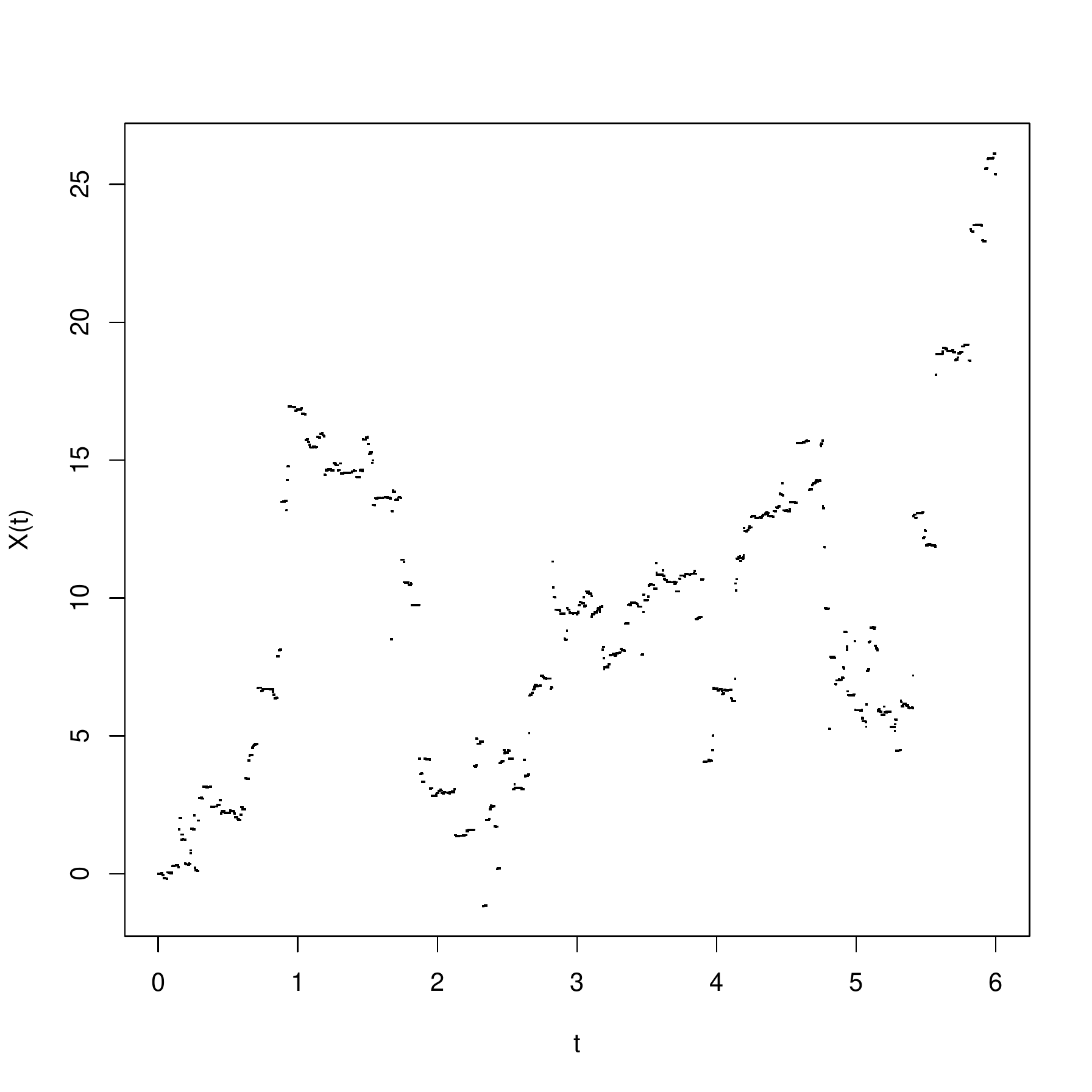}
\caption{$\alpha=0.5$, $f=1$, $\sigma(1)=\sigma(-1)=1$.}
\includegraphics[scale=.37]{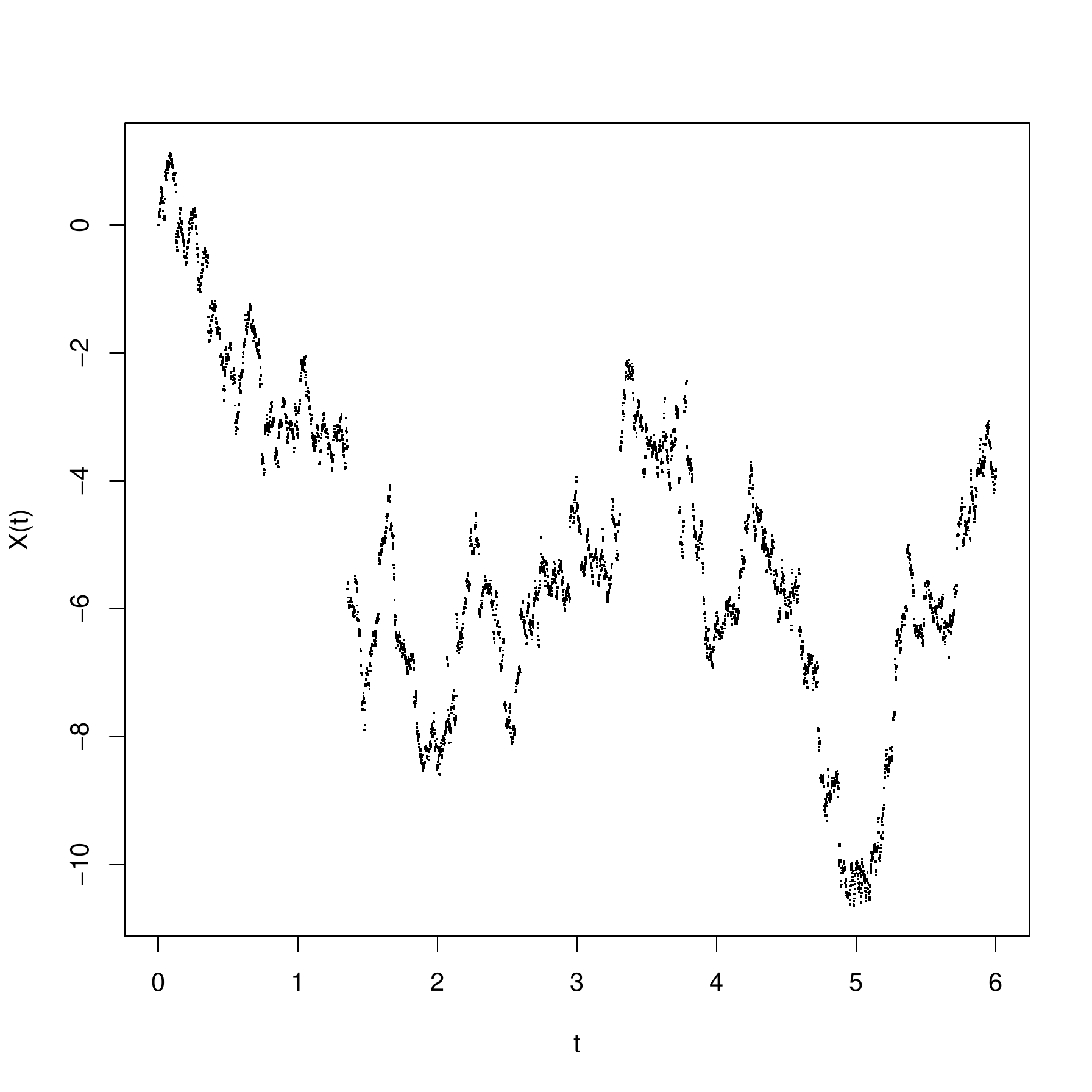}
\caption{$\alpha=1.5$, $f=1$, $\sigma(1)=\sigma(-1)=1$.}
\includegraphics[scale=.37]{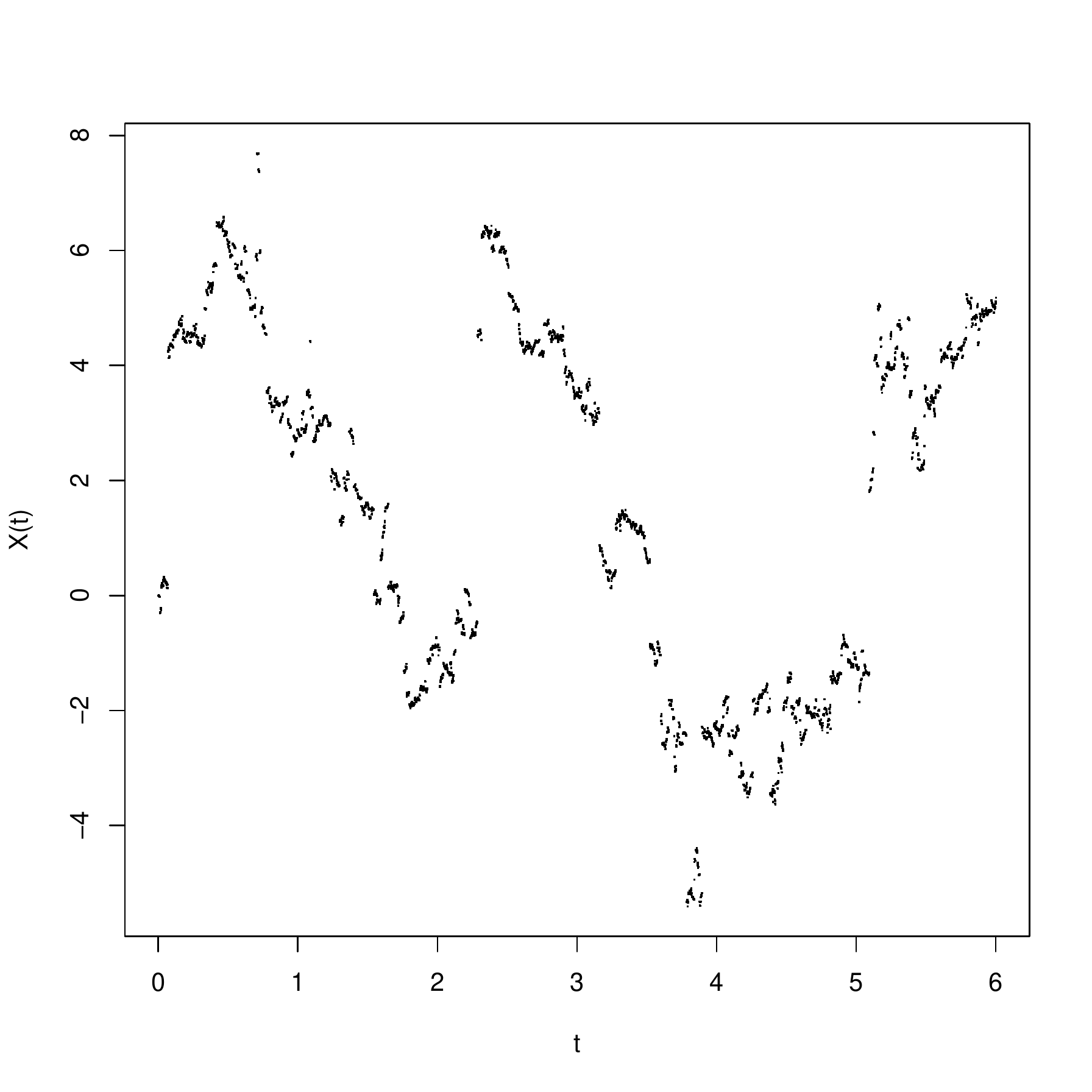}
\caption{$\alpha=1$, $f=1$, $\sigma(1)=\sigma(-1)=1$.}
\end{figure}
\begin{figure}
\centering
\includegraphics[scale=.37]{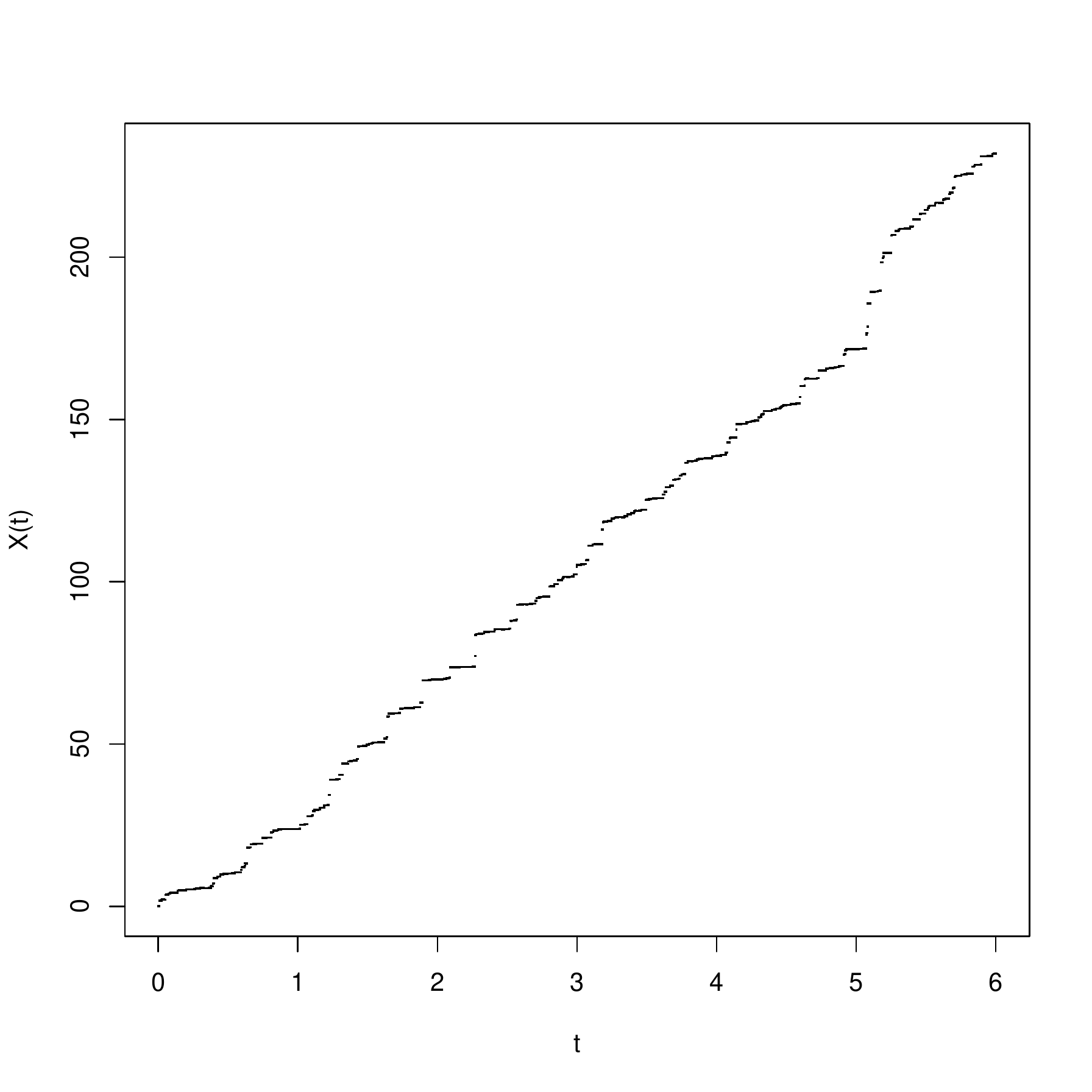}
\caption{$\alpha=0.5$, $f=1$, $\sigma(1)=1$, $\sigma(-1)=0$.}
\includegraphics[scale=.37]{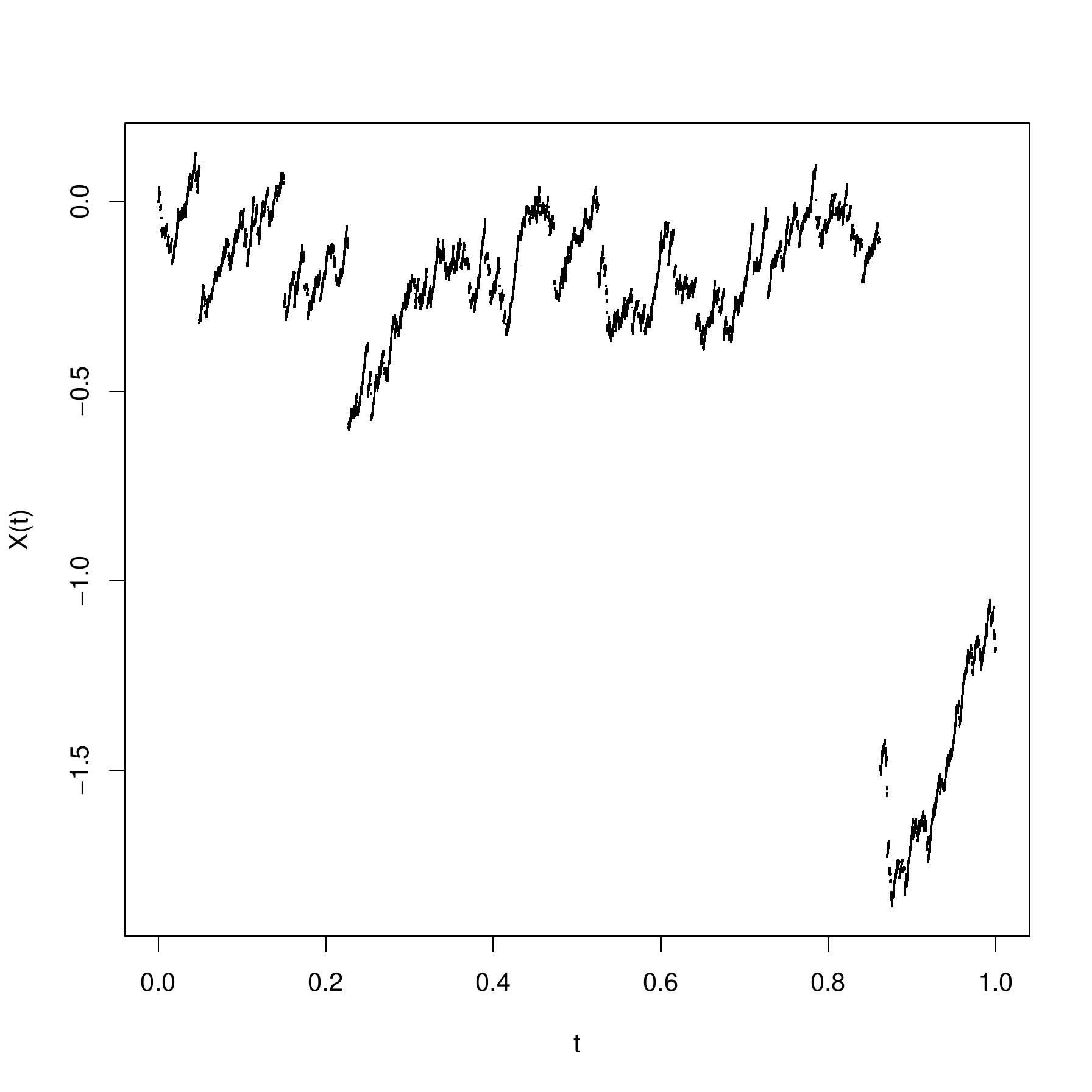}
\caption{$\alpha=1.5$, $f=1$, $\sigma(1)=0$, $\sigma(-1)=1$.}
\includegraphics[scale=.37]{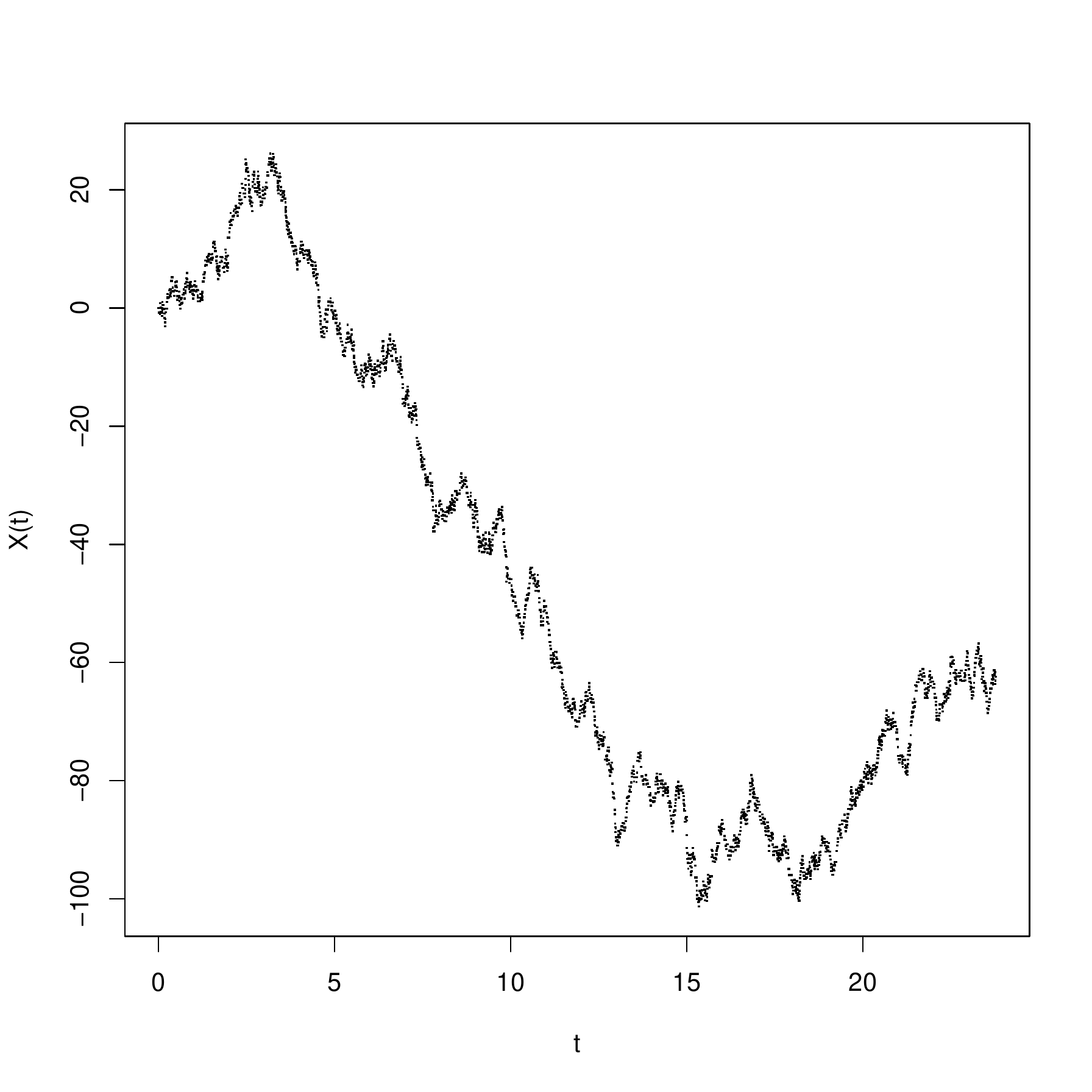}
\caption{$\alpha=1.9$, $f=1$, $\sigma(1)=1$, $\sigma(-1)=1$.}
\end{figure}
\\
\section{Associated processes and examples.}
Here, we are interested in study some related processes
to Lamperti stable distributions (or processes) and give some
examples of Lamperti stable processes which appear in the literature but they are not the main
objects in study. In particular, we study the Ornstein-Uhlenbeck
process and the self-similar additive process related to a
Lamperti stable distribution in the case when the latter is
self-decomposable. We also investigate the parent process of a
Lamperti stable subordinator.
\subsection{Ornstein-Uhlenbeck type processes and self-simlar additive processes.}
Ornstein-Uhlenbeck type processes  appear in  many areas of science, for instance in physics, biology and mathematical finance.  One of the particularity of  these processes  is that its limiting distribution is self-decomposable.  Recall that a random variable $Y$ on $\R^d$, distributed as a Lamperti stable law with characteristics $(\alpha,f, \sigma)$ is self-decomposable if and only if $f\leq \alpha+1/2$.  Therefore, according to Wolfe \cite{wo} and Jurek and Vervaat \cite{jv}, there exists a L\'evy process $Z=(Z_t, t\ge 0)$ on $\R^d$, with $\mathbf{E}_0(\log^+|Z_1|)<\infty$ such that
\[
Y\overset{law}=I:=\int_0^{\infty} e^{-cs}\ud Z_s,
\]
where $c>0$. Consequently, one can define an Ornstein-Uhlenbeck type process driven by $Z$ with initial state $U_0$ and parameter $c>0$, that is the solution of
\[
U_t=U_0+Z_t-c\int_0^t U_s\ud s,
\]
and such that the law of $U_t$ converge towards the law of $Y$ as $t$ goes to $\infty$.  From Theorem 17.5 in \cite{sa}, we have that the process $Z$ has no Gaussian component, its L\'evy measure is given by
\[
\Pi_{Z}(B)=-c\int_{S^{d-1}}\sigma(\ud \xi)\int_{0}^{\infty}\ind_{B}(r\xi)h(r, \xi)\ud r, \qquad B\in \mathcal{B}(\R^d),
\]
where
\[
h(r)=\frac{e^{rf(\xi)}}{(e^r-1)^{\alpha+2}}\Big(re^r(f(\xi)-\alpha-1)+e^r-rf(\xi)-1\Big);
\]
and  linear term
\[
\upsilon=c\eta-\int_{\{\|x\|\ge 1\}}\frac{x}{\|x\|}\Pi_Z(\ud x).
\]
In the one dimensional case, the form of the L\'evy measure of $Z$ is reduced to
\[
\begin{split}
\Pi_Z(\ud x)&=c\bigg(c_+\frac{e^{\beta x}}{(e^x-1)^{\alpha+2}}\Big(x\beta+1-e^x+xe^x(\alpha+1-\beta)\Big)\ind_{\{x>0\}}\\
&-c_-\frac{e^{-\rho x}}{(e^{-x}-1)^{\alpha+2}}\Big(x\rho-1+e^{-x}+xe^{-x}(\alpha+1-\rho)\Big)\ind_{\{x<0\}}\bigg)\ud x,
\end{split}
\]
where $\beta=f(1)$, $\rho=f(-1)$, $c_+=\sigma(\{1\})$ and $c_-=\sigma(\{-1\})$, as usual. In this case, we have another process which is related to  the Lamperti stable distribution $Y$, to $Z$ and to the Ornstein-Uhlenbeck type process $U$. To this end,  recall that in Theorem 16.1 of \cite{sa}, it is showed that  a distribution  is self-decomposable if and only if for any fixed $H>0$, it is the distribution of $V_1$ for some  additive process $V=(V_t,t\geq 0)$ which is self-similar, meaning that for each $k>0$
\[
(V_{kt}, t\geq 0)\overset{d}=\Big(k^HV_t, t\geq 0\Big).
\]
We remark that self-similar additive processes can be used to model space-time scaling random phenomena that can be observed in many areas of science. In particular, they are recently used to  model asset prices and the risk-neutral process  (see for instance \cite{cgmy}) in financial mathematics.

Assume that $V$ is the self-similar additive process associated to $I$, in which case $V_1$ has the same law $I$, and that $H=c$. From Theorem 1 in \cite{jpy}, there are two independent copies of  $Z$ denoted by $Z^{(-)}=(Z^{(-)}_t,t\ge 0)$ and $Z^{(+)}=(Z^{(+)}_t,t\ge 0)$ which are defined by
\[
Z^{(-)}_t\eqdef\int_{e^{-t}}^{1}\frac{\ud V_r}{r^{\gamma}}\qquad\textrm{and} \qquad Z^{(+)}_t\eqdef\int_{1}^{e^{t}}\frac{\ud V_r}{r^{\gamma}}.
\]
The process $V$ can be recovered by
\[
V_r=\left\{
\begin{array}{ll}
{\displaystyle \int_{\log(1/r)}^{\infty} e^{-ct}\ud Z^{(-)}_t} & \textrm{ if }\, 0\leq r\leq 1,\\
{\displaystyle Y+\int_0^{\log (r)} e^{ct}\ud Z^{(+)}_t} & \textrm{ if }\, r\geq 1,
\end{array}
\right.
\]
and moreover $(U^{(+)}_t=e^{-tc}V_{e^t},t\geq 0 )$ is the Ornstein-Uhlenbeck process driven by $Z^{(+)}$ with initial state $I$ and parameter $c$; and $(U^{(-)}_t=e^{tc}V_{e^{-t}},t\geq 0 )$ is the Ornstein-Uhlenbeck process driven by $-Z^{(-)}$ with initial state $I$ and parameter $-c$.
\subsection{Parent process}
Motivated in generating new examples of scale functions, Kyprianou and Rivero \cite{kr} constructed  L\'evy processes with no positive jumps around a given possibly killed  subordinator which plays the role of the descending ladder height process.  One of our aims  is to determine the characteristics of the L\'evy process with no positive jumps whose descending ladder height process is a Lamperti stable subordinator.

Let $X^{L}$ be a Lamperti stable subordinator with characteristics $(\alpha, \beta, \sigma,\theta)$ with zero drift and no killing rate. Since the density of its L\'evy measure is decreasing, then according to Theorem 1 in \cite{kr}, there is $X^{PL}=(X^{PL}_t, t\geq 0)$, a L\'evy process with no positive jumps that we call  {\it the parent process of } $X^L$ whose Laplace exponent is
given by
\[
\psi_{PL}(\lambda)=\lambda\Phi_L(\lambda), \qquad \textrm{for}\quad\lambda\geq 0,
\]
where $\Phi_L$ is the Laplace exponent of $X^L$. Moreover, the process  $X^{PL}$ has no Gaussian coefficient, its L\'evy measure is given by
\begin{equation}
\Pi_{PL}(\ud x)=c_+\frac{e^{-\beta x}}{(e^{-x}-1)^{\alpha+2}}\Big((\alpha+1-\beta)e^{-x}+\beta\Big)\ud x\qquad\textrm{for }\quad x<0,\notag
\end{equation}
and with  linear term
\[
b=\int_{(-\infty,1)}x\Pi_{PL}(\ud x).
\]
Note that $X^{PL}$ oscillates or drifts to $\infty$ according to whether $\Phi_{L}(0)$ is equal zero or strictly positive.  From the form of its L\'evy measure, we deduce that $X^{PL}$ is the sum of two Lamperti stable processes with no positive jumps $X^1$ and $X^{2}$ with characteristics $(\alpha+1, \beta+1, \sigma_1, b_1)$ and $(\alpha+1, \beta, \sigma_2, b_2)$, where $\sigma_{1}(\{1\})=\sigma_{2}(\{1\})=0$,
\begin{align*}
\sigma_{1}(\{-1\})&=c_+(\alpha+1-\beta),& \sigma_{2}(\{-1\})&=c_+\beta,\\
b_1&=\int_{(-\infty,-1)}x\Pi_1(\ud x)-\tilde{a}_{\beta+1},& b_1&=\int_{(-\infty,-1)}x\Pi_2(\ud x)-\tilde{a}_{\beta},
\end{align*}
and $\Pi_1$ and $\Pi_2$ are the respective L\'evy measures of $X^1$ and $X^2$. On the other hand,  the binomial expansion give us
\begin{equation}\label{calpp}
\int_x^\infty \frac{e^{\beta x}}{(e^x-1)^{\alpha+1}}\ud x=e^{-x(\alpha+1-\beta)}\sum_{n=0}^{\infty}\frac{(\alpha+1)_n(\alpha-\beta)_n}{n!(\alpha+1-\beta)_n}e^{-nx},
\end{equation}
which is clearly log-convex on $(0,\infty)$ since it is completely monotone. Hence according to Theorem in 2 \cite{kr}, there is a subordinator $X^{*,L}$ with Laplace exponent $\Phi^*_L$ such that
\[
\Phi_{L}(\lambda)=\frac{\lambda}{\Phi^*_L(\lambda)} \qquad\textrm{for}\quad \lambda\ge 0.
\]
Moreover the subordinator $X^{*,L}$ has no drift and no killing term and the scale function of the parent process $X^{PL}$ is determined
by
\[
W(x)=\int_0^x\Pi^*_L(y,\infty)\ud y,
\]
where $\Pi^*_L$ is the L\'evy measure of $X^{*,L}$. Note that for  $\beta=1$, we have that
\[
\Pi^*_L(y,\infty)=\frac{\alpha}{\Gamma(\alpha)\Gamma(1-\alpha)}(1-e^{-y})^{\alpha-1},\qquad \textrm{for }\quad y>0,
\]
but for $\beta\neq 1$  we do not have an explicit form for $\Pi^*_L$.

The Example 2 in \cite{kr} is related to the Lamperti stable subordinators considered above but with a given killing rate. Let us explain in detail such example in terms of Lamperti stable processes. Take $X^{L,K}$ to be a Lamperti stable subordinator with characteristics $(\alpha, \beta, \sigma,\theta)$ with zero drift and killing rate given by
\[
K=\frac{c_+\Gamma(-\alpha)\Gamma(1-\beta+\alpha)}{\Gamma(1-\beta)}.
\]
According to Kyprianou and Rivero \cite{kr}, there is a subordinator, here now denoted by $Y$ with no drift, no killing rate and L\'evy measure given by
\[
\Pi_Y(\ud x)=\frac{1}{c_+\Gamma^2(1-\alpha)}\bigg((\alpha-\beta)e^{-(\alpha-\beta)x}(e^x-1)^{\alpha-1}+(2-\alpha)\frac{e^{-(\alpha-1-\beta)x}}{(e^x-1)^{2-\alpha}}\bigg)\ud x,
\]
which is the sum of two subordinators, one of which is a Lamperti stable with characteristics $(1-\alpha,\beta+1-\alpha, \sigma_Y)$, where
\[
\sigma_Y(\{1\})=\frac{(2-\alpha)}{c_+\Gamma^2(1-\alpha)},\quad \textrm{and }\quad\sigma_Y(\{-1\})=0.
\]
Moreover, the Laplace exponent of the subordinator $Y$ satisfies that
\[
\phi_Y(\lambda)=\frac{\lambda}{\phi_L(\lambda)},\qquad\textrm{for}\quad \lambda\geq 0.
\]
From the form of $\Pi_Y$, we have the restriction that $\beta<1$.
Thus, his parent process $Y^P$, a spectrally negative L\'evy
process, has Laplace exponent
\[
\psi_{Y^P}(\lambda)=\frac{\lambda^2\Gamma(1-\beta+\lambda)}{\Gamma(1-\beta+\lambda+\alpha)},
\]
which has no Gaussian component and its L\'evy measure satisfies
\[
\Pi_{Y^P}(-\infty,y)=\Pi_Y(\ud y)/\ud y.
\]
According to Kyprianou and Rivero and by (\ref{calpp}), its associated scale function is given by
\[
W_{Y^P}(x)=-Kx+c_+ \sum_{n=0}^{\infty}\frac{(\alpha+1)_n(\alpha-\beta)_n}{n!(\alpha+2-\beta)_n}\Big(1-e^{-(\alpha+2-\beta+n)x}\Big), \quad x\geq 0.
\]
Now $Y^{*,P}$, the parent process of the Lamperti subordinator $X^{L,K}$ with killing rate $K$,  is a spectrally Levy process which drifts to $\infty$, with Laplace exponent
\[
\psi_{Y^{P,*}}(\lambda)=\frac{c_+\Gamma(-\alpha)\lambda\Gamma(\lambda +1 -\beta+\alpha)}{\Gamma(\lambda+1-\beta)},
\]
which has no Gaussian coefficient and whose L\'evy measure satisfies
\[
\Pi_{Y^{P,*}}(\ud x)=\Pi_{PL}(\ud x),\qquad x<0,
\]
with linear term
\[
b=\int_{(-\infty,1)}x\Pi_{PL}(\ud x)-K,
\]
and the associated scale function is given by
\[
W^*(x)=\frac{1}{c_+\Gamma^2(1-\alpha)}\int_0^x e^{-(\alpha-\beta)y}(e^y-1)^{\alpha-1}\ud y.
\]
As the process $X^{PL}$, $Y^{*,P}$ may be seen as the sum of two Lamperti stable process with no positive jumps $Y^1$ and $Y^2$ with characteristics $(\alpha+1, \beta+1, \sigma_1, b_1-K)$ and $(\alpha+1, \beta, \sigma_2, b_2)$, where $\sigma_1,\sigma_2, b_1$ and $b_2$ are defined as above.

It is important to note that the above example have been recently used for the risk neutral stock price model by Eberlein and Madan \cite{em}.
\subsection{Examples.}
Examples of Lamperti stable processes  appear in the literature at least in the papers mentioned  in the introduction (\cite{cc, ckp, kp, pp}) but they also appear (in a hidden
way) in many other recent works. We will give a quick overview of
some of them, not pretending to be exhaustive in this list.

In \cite{by}, we find two examples related to the factorization
$$\mathbf{e}\overset{law}=\mathbf{e}^{\alpha}\tau_{\alpha}^{-\alpha}.$$
where ${\mathbf{e}}$ is an exponential variable independent of the
$\alpha$-stable variable $\tau_{\alpha}$. The first of them is related with the exponential functional of a killed subordinator $Z^1$ whose  Laplace
exponent is given by
$$\phi_1 (\lambda) = \frac{\Gamma(\alpha \lambda+1)}{\Gamma(\alpha(\lambda-1)+1)}.$$ It is easy to see that it is related to the  Laplace exponent $\Phi_L$ of a Lamperti stable subordinator $X^{L}$ with characteristics $(\alpha,\alpha,\sigma,\theta)$, zero drift, and $\sigma(\{1\})=\alpha/\Gamma(1-\alpha)$.
The relationship between both Laplace exponents is
$$\phi_1(\lambda)= \Phi_L(\alpha\lambda)
+\frac{1}{\Gamma(1-\alpha)}.$$ This subordinator is also studied in Rivero
\cite{riv}, where the author finds its renewal density and other related
computations.

The Laplace exponent of the second subordinator, here denoted by $Z^2$, is given by
$$\phi_2(\lambda)=\lambda\frac{\Gamma(\alpha(\lambda-1)+1)}{\Gamma(\alpha \lambda
+1)},$$ and can be expressed in terms of the Laplace exponent $\Phi_{L,2}$
of a Lamperti stable subordinator $X^{L,2}$ with characteristics
$(1-\alpha,1,\sigma,\theta)$, and zero drift where $\sigma(\{1\})=\alpha/\Gamma(1-\alpha)$.
The relation between them is
$$\phi_2(\lambda)= \alpha\Phi_{L,2}(\alpha\lambda).$$

In both cases  this allows us to compute the law of the exponential
functional of  $\alpha X^L$ and $\alpha X^{L,2}$ in terms of the one of
$Z^1$ and $Z^2$, respectively.

There is another example in \cite{by} which is related to the factorization

$$\mathbf{e}\overset{law}=\gamma_s^{\alpha}J_s^{(\gamma)},$$
where $s\ge \alpha$, $\gamma_s$ is a Gamma r.v. with parameter $s$ and $J_s^{(\gamma)}$ denotes a certain r.v which is independent of $\gamma_s$. In this case, the killed subordinator related to the exponential functional which has the same moments as the  $\gamma_s$, can be expressed as the
sum of two independent Lamperti stable procesess. In \cite{riv}
further calculation are carried over concerning this subordinator.

 In the paper \cite{svy} in section 5.3, the authors found the L\'evy measure
of the inverse of the local time at 0 of an Ornstein Uhlenbeck process driven by a standard Brownian motion and parameter $\gamma>0$.
This measure is
\begin{equation}
\nu(t)=\frac{\gamma^{3/2}e^{\gamma t/2}}{\sqrt{2\pi}(\sinh(\gamma
t))^{3/2}}=\frac{(2\gamma)^{3/2}e^{2\gamma
t}}{\sqrt{2\pi}(e^{2\gamma t}-1)^{3/2}}\notag
\end{equation}
and the corresponding Laplace exponent is computed. It is related to a
Lamperti stable distribution with characteristics
$(1/2,1,\sqrt{\gamma/\pi}).$

This computation as well as the three former examples can be carried
out by recognizing that behind those measures there is a related
Lamperti stable distribution and applying our Theorem 1 to calculate
the corresponding Laplace exponent.

In the papers \cite{ckp} , \cite{kp}, \cite{pp} the main processes in study
are Lamperti stable processes. All these papers share the property
that many useful explicit calculations are be carried out. This is,
we believe, the main advantage of this class: being at the same time
a good model for many situations, allowing simulation of the paths
as well as many explict calculation to be carried on.\\
\\ \

\noindent {\bf Acknowledgements.} This research was supported by EPSRC grant EP/D045460/1, CONACYT
grant (200419), and the project PAPIITT-IN120605. We are much
indebted to Andreas Kyprianou, Victor Rivero, and Ger\'{o}nimo Uribe for many fruitful
discussions on L\'evy processes with no positive jumps and parent
processes.

\end{document}